\RequirePackage[leqno]{amsmath}
\documentclass[leqno]{amsart}
\usepackage[utf8]{inputenc}     % Accept diferent input encodings
\usepackage{ifthen}  
\usepackage{multicol}
\usepackage{color, graphics}
\usepackage{pgfplots}
\usepackage{graphicx, multicol, latexsym, amsmath, amssymb}
\usepackage{subcaption}
\usepackage{amsfonts}
\usepackage{enumitem}
\setlist[itemize]{topsep=0pt,after=\vspace{1.5\baselineskip}}
 \usepackage[colorlinks=true]{hyperref}
\usepackage{empheq}
\usepackage[T1]{fontenc}
\usepackage{tabularx}
\usepackage[leftcaption]{sidecap}
\sidecaptionvpos{figure}{m}
\usepackage{tikz}
\usepackage[margin=0.5in]{geometry}

\setlist[itemize]{noitemsep, topsep=0pt}

\def\R{\mathbb R} \def\N{\mathbb N}

\def\R{\mathbb R} \def\N{\mathbb N} 
\def\TM{T_{max}} 
\def
\@cite
#1#2{[{{\bfseries #1}\if@tempswa , #2\fi}]}
\newtheorem{theorem}{Theorem}[section]

\newtheorem{lemma}[theorem]{Lemma}

\newtheorem{remark}{Remark}

\RequirePackage[normalem]{ulem} %DIF PREAMBLE
\RequirePackage{color}\definecolor{RED}{rgb}{1,0,0}\definecolor{BLUE}{rgb}{0,0,1} %DIF PREAMBLE
\providecommand{\DIFaddbegin}{} %DIF PREAMBLE
\providecommand{\DIFaddend}{} %DIF PREAMBLE
\providecommand{\DIFdelbegin}{} %DIF PREAMBLE
\providecommand{\DIFdelend}{} %DIF PREAMBLE
\providecommand{\DIFaddbeginFL}{} %DIF PREAMBLE
\providecommand{\DIFaddendFL}{} %DIF PREAMBLE
\providecommand{\DIFdelbeginFL}{} %DIF PREAMBLE
\providecommand{\DIFdelendFL}{} %DIF PREAMBLE
\newcommand{\DIFscaledelfig}{0.5}
\RequirePackage{settobox} %DIF PREAMBLE
\RequirePackage{letltxmacro} %DIF PREAMBLE
\newsavebox{\DIFdelgraphicsbox} %DIF PREAMBLE
\newlength{\DIFdelgraphicswidth} %DIF PREAMBLE
\newlength{\DIFdelgraphicsheight} %DIF PREAMBLE
\LetLtxMacro{\DIFOincludegraphics}{\includegraphics} %DIF PREAMBLE
\newcommand{\DIFaddincludegraphics}[2][]{{\color{blue}\fbox{\DIFOincludegraphics[#1]{#2}}}} %DIF PREAMBLE
\newcommand{\DIFdelincludegraphics}[2][]{% %DIF PREAMBLE
\sbox{\DIFdelgraphicsbox}{\DIFOincludegraphics[#1]{#2}}% %DIF PREAMBLE
\settoboxwidth{\DIFdelgraphicswidth}{\DIFdelgraphicsbox} %DIF PREAMBLE
\settoboxtotalheight{\DIFdelgraphicsheight}{\DIFdelgraphicsbox} %DIF PREAMBLE
\scalebox{\DIFscaledelfig}{% %DIF PREAMBLE
\parbox[b]{\DIFdelgraphicswidth}{\usebox{\DIFdelgraphicsbox}\\[-\baselineskip] \rule{\DIFdelgraphicswidth}{0em}}\llap{\resizebox{\DIFdelgraphicswidth}{\DIFdelgraphicsheight}{% %DIF PREAMBLE
\setlength{\unitlength}{\DIFdelgraphicswidth}% %DIF PREAMBLE
\begin{picture}(1,1)% %DIF PREAMBLE
\thicklines\linethickness{2pt} %DIF PREAMBLE
{\color[rgb]{1,0,0}\put(0,0){\framebox(1,1){}}}% %DIF PREAMBLE
{\color[rgb]{1,0,0}\put(0,0){\line( 1,1){1}}}% %DIF PREAMBLE
{\color[rgb]{1,0,0}\put(0,1){\line(1,-1){1}}}% %DIF PREAMBLE
\end{picture}% %DIF PREAMBLE
}\hspace*{3pt}}} %DIF PREAMBLE
} %DIF PREAMBLE
\LetLtxMacro{\DIFOaddbegin}{\DIFaddbegin} %DIF PREAMBLE
\LetLtxMacro{\DIFOaddend}{\DIFaddend} %DIF PREAMBLE
\LetLtxMacro{\DIFOdelbegin}{\DIFdelbegin} %DIF PREAMBLE
\LetLtxMacro{\DIFOdelend}{\DIFdelend} %DIF PREAMBLE
\DeclareRobustCommand{\DIFaddbegin}{\DIFOaddbegin \let\includegraphics\DIFaddincludegraphics} %DIF PREAMBLE
\DeclareRobustCommand{\DIFaddend}{\DIFOaddend \let\includegraphics\DIFOincludegraphics} %DIF PREAMBLE
\DeclareRobustCommand{\DIFdelbegin}{\DIFOdelbegin \let\includegraphics\DIFdelincludegraphics} %DIF PREAMBLE
\DeclareRobustCommand{\DIFdelend}{\DIFOaddend \let\includegraphics\DIFOincludegraphics} %DIF PREAMBLE
\LetLtxMacro{\DIFOaddbeginFL}{\DIFaddbeginFL} %DIF PREAMBLE
\LetLtxMacro{\DIFOaddendFL}{\DIFaddendFL} %DIF PREAMBLE
\LetLtxMacro{\DIFOdelbeginFL}{\DIFdelbeginFL} %DIF PREAMBLE
\LetLtxMacro{\DIFOdelendFL}{\DIFdelendFL} %DIF PREAMBLE
\DeclareRobustCommand{\DIFaddbeginFL}{\DIFOaddbeginFL \let\includegraphics\DIFaddincludegraphics} %DIF PREAMBLE
\DeclareRobustCommand{\DIFaddendFL}{\DIFOaddendFL \let\includegraphics\DIFOincludegraphics} %DIF PREAMBLE
\DeclareRobustCommand{\DIFdelbeginFL}{\DIFOdelbeginFL \let\includegraphics\DIFdelincludegraphics} %DIF PREAMBLE
\DeclareRobustCommand{\DIFdelendFL}{\DIFOaddendFL \let\includegraphics\DIFOincludegraphics} %DIF PREAMBLE
\title[Boundedness in a Keller--Segel model with nonlinear production] %Use the shortened version of the full title
      {%A nonlocal reaction chemotaxis model with arbitrarily positive diffusion exponent
Boundedness for a fully parabolic Keller--Segel model with sublinear segregation and superlinear aggregation}
\author[S. Frassu and G. Viglialoro]{}
\keywords{Chemotaxis, 
	Chemotactic sensitivity, Global existence, Nonlinear production. \\
	\textit{$^*$Corresponding author}: giuseppe.viglialoro@unica.it}
\begin{document}
\maketitle
%\tableofcontents

\centerline{\scshape Silvia Frassu \and Giuseppe Viglialoro$^{*}$}
\medskip
{
% \footnotesize
% \centerline{$^1$Institut für Mathematik}
% \centerline{Universität  Paderborn }
% \centerline{Warburge
%r  Str.    100,  33098  Paderborn (Germany)}
% \centerline{$^\natural$School of Control Science and Engineering}
% \centerline{Shandong University}
% \centerline{Jinan, Shandong, 250061. (P. R. China)}
 \medskip
 \centerline{Dipartimento di Matematica e Informatica}
 \centerline{Universit\`{a} di Cagliari}
 \centerline{Via Ospedale 72, 09124. Cagliari (Italy)}
%  \centerline{E-mail: giuseppe.viglialoro@unica.it}
 \medskip
}
\bigskip
\begin{abstract}
This work deals with a fully parabolic chemotaxis model with nonlinear production and chemoattractant. The problem is formulated on a bounded domain and, depending on a specific interplay between the coefficients associated to such production and chemoattractant,  we establish that the related initial-boundary value problem has a unique classical solution which is uniformly bounded in time. To be precise, we study this zero-flux problem
\begin{equation}\label{problem_abstract}
\tag{$\Diamond$}
\begin{cases}
u_t= \Delta u - \nabla \cdot (f(u) \nabla v) & \text{ in } \Omega \times (0,T_{max}),\\
v_t=\Delta v-v+g(u)  & \text{ in } \Omega \times (0,T_{max}),\\
%u(x,0)=u_0(x), \; v(x,0)=v_0(x) & x \in \bar\Omega,
\end{cases}
\end{equation}
where $\Omega$ is a bounded and smooth domain of $\R^n$, for $n\geq 2$, and $f(u)$ and $g(u)$ are reasonably regular functions generalizing, respectively, the prototypes $f(u)=u^\alpha$ and $g(u)=u^l$, with proper $\alpha, l>0$. After having shown that  any sufficiently smooth $ u(x,0)=u_0(x)\geq 0, \, v(x,0)=v_0(x)\geq 0$ emanate a unique classical and nonnegative solution $(u,v)$ to problem \eqref{problem_abstract}, which is defined on $\Omega \times (0,T_{max})$ with $T_{max}$  
 denoting the maximum time of existence, we establish that for any $l\in (0,\frac{2}{n})$ and $\frac{2}{n}\leq \alpha<1+\frac{1}{n}-\frac{l}{2}$, $T_{max}=\infty$ and $u$ and $v$ are actually  uniformly bounded in time.

This paper is  in line with the contribution by Horstmann and Winkler (\cite{HorstWink}) and, moreover, extends the result by Liu and Tao (\cite{LiuTaoFullyParNonlinearProd}).  Indeed, in the first work it is proved that for $g(u)=u$ the value $\alpha=\frac{2}{n}$ represents the critical blow-up exponent to the model, whereas in the second, for $f(u)=u$, corresponding to $\alpha=1$, boundedness of solutions is shown under the assumption $0<l<\frac{2}{n}.$  
\end{abstract}
\begingroup
\hypersetup{hidelinks}
\tableofcontents
\endgroup
\section{Introduction and motivations}\label{IntroSection}
Most of this article is dedicated to the following Cauchy boundary problem
\begin{equation}\label{problem}
\begin{cases}
u_t= \Delta u - \nabla \cdot (f(u) \nabla v) & \text{ in } \Omega \times (0,T_{max}),\\
v_t=\Delta v-v+g(u)  & \text{ in } \Omega \times (0,T_{max}),\\
u_{\nu}=v_{\nu}=0 & \text{ on } \partial \Omega \times (0,T_{max}),\\
u(x,0)=u_0(x), \; v(x,0)=v_0(x) & x \in \bar\Omega,
\end{cases}
\end{equation}
defined  in a bounded and smooth domain $\Omega$ of $\R^n$, with $n\geq 2$, and formulated through some  functions $f=f(s)$ and $g=g(s)$, sufficiently regular in their argument $s\geq 0$, and further regular initial data  $u_0(x)\geq 0$ and $v_0(x)\geq 0.$  Additionally,  the subscription $(\cdot)_\nu$ indicates the outward normal derivative on $\partial \Omega$ whereas $T_{max}$ the maximum time up to which solutions to the system are defined. 

The two partial differential equations appearing above generalize  
\begin{equation}\label{problemOriginalKS}
\begin{cases}
u_t= \Delta u - \nabla \cdot (u \nabla v) & \text{ in } \Omega \times (0,T_{max}),\\
v_t=\Delta v-v+u  & \text{ in } \Omega \times (0,T_{max}),
\end{cases}
\end{equation}
proposed in the pioneer papers by Keller and Segel (\cite{K-S-1970,Keller-1971-MC}) to model the dynamics of populations (as for instance cells or bacteria) arising in mathematical biology. Precisely, by indicating with  $u = u(x, t)$ a certain particle density at the position $x$ and at the time $t$, the equations describe how the aggregation impact from the coupled cross term $u \nabla v$, related to the chemical signal $v=v(x,t)$ (initially distributed accordingly to the law $v_0(x)=v(x,0)$, as in \eqref{problem}), may contrast the natural diffusion  (associated to the Laplacian operator, $\Delta u$) of the cells, organized at the initial time through the configuration $u_0(x)=u(x,0)$. In particular, such an attractive impact might influence the motion of the cells  so strongly even to lead the  system to its chemotactic collapse (blow-up at finite time with appearance of $\delta$-formations for the particle density). In the literature there are many contributions dedicated to the comprehension of this phenomena. In this regard, in  \cite{JaLu,Nagai} the reader can find an extensive theory dealing with the existence and properties of global, uniformly bounded or blow-up (local) solutions to the Cauchy problem associated to \eqref{problemOriginalKS}, and endowed with homogeneous Neumann boundary conditions (exactly as in \eqref{problem}, and biologically modeling an impermeable domain), especially in terms of the initial mass of the particle distribution, i.e.,  $m=\int_\Omega u_0(x)dx.$ Indeed, the mass of the bacteria, preserved in time for this model, appears as a critical parameter (see, for instance,  \cite{HerreroVelazquez,Nagai,WinklAggre}); more exactly, for $n\geq 2$,  the value $m_c=4\pi$ establishes that when the diffusion overcomes the self-attraction ($m<m_c$), global in time solutions are expected, whereas when the self-attraction dominates the diffusion ($m>m_c$), blow-up solutions at finite time may be detected. 

The size of the initial distribution is not the only factor capable to influence the chemotactic behavior of the cells toward their self-organization. Other elements may also take sensitively part in this process; the impacts of the diffusion and/or of the chemoattractant, weaker or stronger (for instance, if in \eqref{problemOriginalKS} the cross-diffusion term $u\nabla v$ is replaced by $\chi u \nabla v$, for some $\chi>0$, then even for initial distribution $u_0$ with subcritical mass, the system exhibits blow-up at finite time whenever $\chi$ increases), the presence of external sources affecting the cells' density, or the law of the signal production from the chemical, dictated by the cells themselves: to a high (low) segregation corresponds a high disorganization (organization) in the motion of the particles. Herein we are interested in the analysis concerning the \textit{mutual interplay between the actions from the chemoattractant and the segregation rates.} In such sense problem \eqref{problem} is an example of chemotaxis model combining these aspects, exactly as specified: the chemosensitivity function $f(u)$ describes how the population aggregates, through the interaction with the chemical, and directs its movement in the direction of the gradient of $v$. In our problem $f(u)$ generalizes $u^\alpha$, for some $\alpha$ even covering superlinear powers;  as said, the larger $\alpha$ the higher is the attraction between each cell, leading the system to undesired instabilities. On the other hand, the second equation indicates that chemical signal is produced according to the law of $g(u)$, which as well has as prototype $u^l$, with $l$ smaller than 1. Naturally, this has a segregation impact on the model weaker than the case with $g(u)=u$, especially at large particle densities (see \cite{Hillen-2009-UGP,Myerscough-1998-PFG,Murray1} an related references therein). As conceivable, the gathering phenomena characterizing the original model is dampened and more smoothness to the system is more efficiently than supplied. 

Before giving our precise objectives in respect of the analysis above developed, let us mention the result which, mainly, inspires and justifies this investigation:  For $g(u)=u$ and $f(u)\cong u^\alpha$, it is shown that the value $\alpha=\frac{2}{n}$ decides whether  model \eqref{problem} manifests or not blow-up scenarios.  Specifically, with $n\geq 2$, for $\alpha\in (0,\frac{2}{n})$ all solutions are global and uniformly bounded, whereas the same does not apply for $\alpha>\frac{2}{n}$. In fact, (a) for $\alpha>2$ and any $n\geq 2$, (b) for  $\alpha \in (1,2)$, $n\in \{2,3\}$ and technical assumptions on $f$, (c) for $\alpha \in (\frac{2}{n},1)$ and $n\in \{2,3\}$ or (d) for $\alpha \in (\frac{2}{n},\frac{2}{n-2})$ and $n\geq 4$ (also in this case combined with further assumptions on $f$),  there are initial data $(u_0,v_0)$ emanating  unbounded solutions (see \cite{HorstWink}) 

By continuing within the confines of Keller--Segel models with linear production, when the diffusion is not linear, i.e. $\Delta u=\nabla \cdot \nabla u$ reads $\nabla \cdot (D(u)\nabla u)$,  for $n\geq 2$ the asymptotic behavior of the ratio $\frac{f(u)}{D(u)}\cong u^\alpha$ for large values of
$u$ indicates that if $\alpha\in (0,\frac{2}{n})$ any $(u_0,v_0)$ produces uniformly bounded classical solutions to problem \eqref{problem} (see \cite{TaoWinkParaPara}), whilst for $\alpha>\frac{2}{n}$ blow-up solutions either in finite of infinite time can be constructed, even for arbitrarily small initial data (see \cite{WinVolumeFill}). Moreover, the insight about the quantitative role of the diffusion of the cells on the evolution of the model reads as follows: for $D(u)\cong u^{m-1}$ and $f(u)\cong u^\alpha$,  $m,\alpha\in \R$, it is established in \cite{CieslakStinnerFiniteTime,CieslakStinnerNewCritical} (for the fully parabolic case) and in \cite{WinDj} (for the simplified parabolic-elliptic one) that $\alpha<m+\frac{2}{n}-1$ is condition sufficient and necessary in order to ensure global existence and boundedness of solutions. (For completeness, we also refer to \cite{MarrasNishinoViglialoro}, where an estimate for the blow-up time of unbounded solutions to the simplified model is derived.) Unlike the case where $D(u)=1$ and $f(u)=u$ where the critical mass $m_c$ is $n$-independent, the above criterion implies that the size of the initial mass may have no crucial role on the existence of global or local-in-time solutions to nonlinear diffusion chemotaxis-systems. Conversely, the key factor is given by some specific interplay between the coefficients $m,\alpha$ and the dimension $n$; in particular, as mentioned, this is especially observed at high dimensions, for which a magnification of the diffusion parameter is required to compensate instability effects. A similar consideration, appropriately reinterpreted in that context, will be given below when the exponent associated to the nonlinear signal production is introduced.    

Complementary, as far as nonlinear segregation chemotaxis models are concerned, when in problem \eqref{problem}  the case $f(u)=u$  is considered, uniformly boundedness of all its solutions is proved in \cite{LiuTaoFullyParNonlinearProd} for $g(u)\cong u^l$, with  $0<l<\frac{2}{n}$. Moreover, by resorting to a simplified parabolic-elliptic version in spatially radial contexts, for $f(u)=u$ and the second equation reduced to $0=\Delta v-\mu(t)+g(u)$, with $g(u)\cong u^l$ and $\mu(t)=\frac{1}{|\Omega|}\int_\Omega g(u(\cdot,t))$, it is known  (see \cite{WinklerNoNLinearanalysisSublinearProduction}) that  the same conclusion on the boundedness continues valid for any $n\geq 1$ and $0<l<\frac{2}{n}$, whereas for $l>\frac{2}{n}$ blow-up phenomena may appear. 
\section{Presentation of the main result and comparison with a simplified model. Plan of the paper}
\subsection{Claim of the main result} 
In accordance to what discussed above, we are interested in addressing situations concerning system \eqref{problem} that, to our knowledge, are not yet discussed in the literature. To this aim, from now these assumptions, respectively identifying the actions associated to the  chemoattractant and to the segregation of the chemical signal, are fixed:
\begin{equation}\label{f}
f \in C^2([0,\infty)), \quad f(0)=0\quad  \textrm{and  }\quad f(s)\leq Ks^{\alpha},\quad \textrm{for some}\quad K,\alpha>0 \quad \textrm{and all } s \geq 0,
\end{equation}
and
\begin{equation}\label{AssumptionOnSegregationG}
g \in C^1([0,\infty))\quad \textrm{and  }\quad 0\leq g(s)\leq K_0 s^l, \quad \textrm{for some}\quad K_0,l>0 \quad \textrm{and all } s \geq 0.
\end{equation}
In particular, with a specific view to what mentioned for \cite{HorstWink}, linear productions of the chemical may be sufficient to emanate blow-up solutions when the impact from the chemoattractat, favoring gatherings in the motion of the species,  is superquadratic, in any dimension,  superlinear and subquadratic, in low dimensions, and sublinear in higher. Thus seems meaningful the following question: 
\begin{center}
\begin{changemargin}{1.5cm}{1.5cm} 
$\circ$ May a sublinear signal segregation of the chemical  enforce globability of solutions  for superlinear chemosensitivitiy even in high dimensions?  
\end{changemargin}
\end{center}
Our result positively addresses this issue in the sense that
 \textit{independently of the initial data, by weakening in an inversely proportional way to the dimension the impact associated to the production rate of the chemical, the uniform-in-time boundedness of solutions to model \eqref{problem} is ensured, even for superlinear thrusts from the chemoattractant.}

What said is formally claimed in this 
\begin{theorem}\label{MainTheorem}   
	%{\textcolor{blue}{
	Let $\Omega$ be a bounded and smooth domain of $ \R^n$, with $n \geq 2$. Moreover, let $f$ and $g$ fulfill \eqref{f} and \eqref{AssumptionOnSegregationG}, respectively, with $l\in (0,\frac{2}{n})$ and $\alpha$ satisfying 
	\begin{equation}\label{AssumptionAlphaElleN}
	%	\begin{cases}
	\frac{2}{n}\leq \alpha <1+\frac{1}{n}-\frac{l}{2}.
	%&  \textrm{if }\, n=2,\\
	%	\frac{2}{n}\leq \alpha \leq 1&  \textrm{if }\, n\geq 3.\\
	%	\end{cases}
	\end{equation}
	Then, for any nontrivial $(u_0,v_0)\in C^0(\bar{\Omega})\times C^{1}(\bar\Omega)$, with $u_0\geq 0$ and $v_0\geq 0$ on $\bar{\Omega}$, there exists a unique pair of nonnegative functions $(u,v)\in (C^0(\bar{\Omega}\times [0,\infty))\cap  C^{2,1}(\bar{\Omega}\times (0,\infty)))^2$ 
	%	$$
	%	\begin{aligned}
	%	u &\in  C^0(\bar{\Omega}\times [0,\infty))\cap  C^{2,1}(\bar{\Omega}\times (0,\infty)),\\
	%	v &\in  C^0(\bar{\Omega}\times [0,\infty))\cap  C^{2,1}(\bar{\Omega}\times (0,\infty)),
	%	\end{aligned}
	%	$$
	which solve problem \eqref{problem} and satisfy for some $C>0$  
	\begin{equation*}
	\lVert u (\cdot,t)\rVert_{L^\infty(\Omega)}+\lVert v (\cdot,t)\rVert_{L^{\infty}(\Omega)}\leq C\quad \textrm{for all }t >0.
	\end{equation*}
\end{theorem}
% , let us first fix these mutual relations on the parameters $\alpha,\beta,\eta,m>0$,  also in term of $l:=\frac{2n}{n-2}$, a precise value related to the space dimension $n$, which herein we suppose to satisfy $n\geq 3$: 
%\begin{align}
%\label{b}\tag{\textsf{\textsf{H1}}}
%\alpha+\beta>\max\{(4 l - 4 - m l)/(p - 2), (2 \eta l - m l - 2 \eta)/(l - 
%    2) , m , 2 - m l/(l - 2)\}.
%\end{align}
%
%\begin{align}
%\label{b1}\tag{\textsf{\textsf{H2}}}
%\alpha+\beta>\max\{(4 l - 4 - m l)/(p - 2), (2 \eta l - m l - 2 \eta)/(l - 
%    2) , 2 - m l/(l - 2)\}.
%\end{align}
%These represent exactly our main assertions concerning the attractive-repulsive case ($\lambda=1$): 
\begin{remark}
We make this considerations:
	\begin{itemize}
\item 		For $\alpha=1$ assumption \eqref{AssumptionAlphaElleN} is simplified into $l\in (0,\frac{2}{n})$. In particular, our analysis is an extension of that developed in \cite{LiuTaoFullyParNonlinearProd}, in the sense that Theorem \ref{MainTheorem} recovers \cite[Theorem 1.1]{LiuTaoFullyParNonlinearProd} when $f(u)=u$ in problem \eqref{problem}. 
\item From $l\in (0,\frac{2}{n})$, the comparison with the limit linear signal production model for system \eqref{problem}  makes sense only in two-dimensional settings; for  $l=1$ the upper bound in assumption \eqref{AssumptionAlphaElleN} reads $\alpha<1$, and Theorem \ref{MainTheorem} is consistent with \cite[Theorem 4.1]{HorstWink}. 
\item Since \cite[Theorem 4.1]{HorstWink} is applicable for any $n\geq 2$ whenever  $\alpha<\frac{2}{n}$ and $l=1$, a fortiori it holds true for $l \in (0,\frac{2}{n})$; this is the sole reason why we consider in our analysis $\alpha\geq \frac{2}{n}$.
\item Considering that for linear production and nonlinear diffusion (with parameter $m$) models we discussed that the condition for boundedness reads $\alpha<m+\frac{2}{n}-1$, from assumption \eqref{AssumptionAlphaElleN}  one can observe that the parameter $l$ associated to the nonlinear segregation plays an opposite role with respect $m$: given $\alpha$, for high values of $n$, smaller (larger) values of $l$ ($m$) are needed to ensure globability and boundedness. 
	\end{itemize}
\end{remark}
	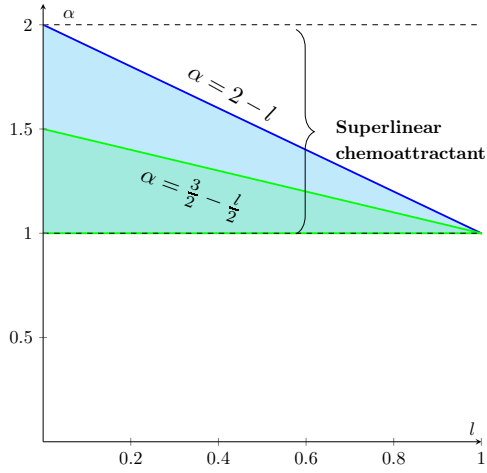
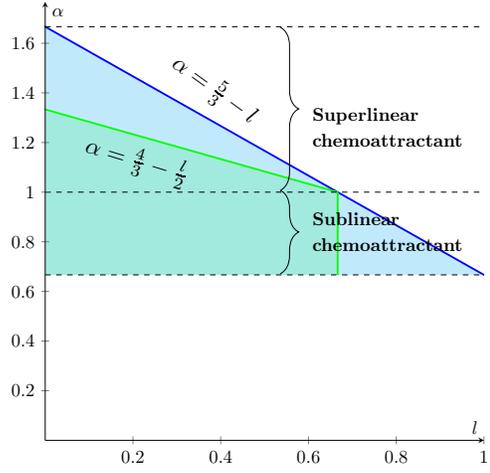
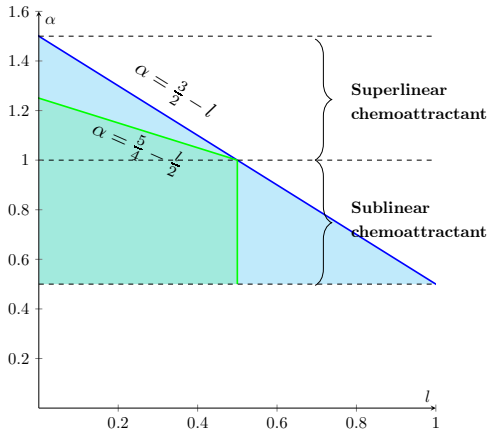
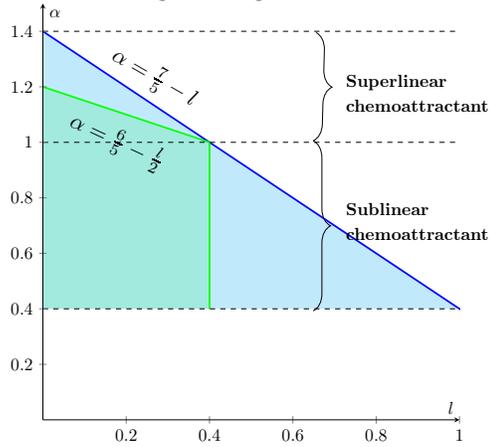
\begin{figure}[h!!]
	\centering 
	%% n=2 and n=3 
	\begin{subfigure}[b]{0.35\linewidth}
		\resizebox{1\textwidth}{!}{
			\begin{tikzpicture}
			\begin{axis}[
			xlabel={$l$},ylabel={$\;\;\, \alpha$},
			xmin=0, xmax=1,
			ymin=0, ymax=1.1+2/2,
			axis y line=center, axis x line=center,
			width=10cm, height=10cm,axis on top]
			\addplot[line width=0.25pt,color=white,fill=green!20,opacity=0.75] coordinates{(0,2/2) (0,1+1/2) (2/2,2-2/2)} 
			node[rotate=-29,right=-5.5cm,above=-0.30cm,scale=1.3,black,opacity=1] {$\alpha=2-l$};
			\addplot[line width=0.25pt,color=white,fill=cyan,opacity=0.25] coordinates{(0,2/2) (0,1+2/2) (1,2/2)} 
			node[rotate=-20,right=-5.5cm,above=-1.6cm,scale=1.3,black,opacity=1] {$\alpha=\frac{3}{2}-\frac{l}{2}$};
			\addplot[line width=1pt,domain=-5:2,samples=2,blue] ({x},{1+2/2-x}); 
			\addplot[line width=1pt,domain=-2:2,samples=2,green] ({x},{1+1/2 -x/2}); 
			\addplot[line width=1pt,domain=-2:2,samples=2,green] ({x},{1}); 
			%Linea alpha=1
			\addplot[line width=0.1pt,dashed,domain=0:1,samples=2,black] ({x},{1}); 	
			\addplot[line width=0.1pt,dashed,domain=0:1,samples=2,black] ({x},{1+2/2}); 		
			\end{axis}
			\filldraw[black] (5.5,6.3) node[anchor=north west] {{\bf{Superlinear}}};
			\filldraw[black] (5.5,5.8) node[anchor=north west] {{\bf{chemoattractant}}};
			%Parentesi
			\draw
			[decorate,decoration={brace,amplitude=10pt,mirror},xshift=-4pt,yshift=0pt]
			(5.,4) -- (5.,7.928)node [black,midway,xshift=13pt] {};
			\end{tikzpicture}
		}
		\caption{Case $n=2$. PP: $l\in (0,1)$, $1\leq \alpha<\frac{3}{2}-\frac{l}{2}$. PE:  $l\in (0,1)$, $1\leq \alpha<2-l$.} \label{fig:}
	\end{subfigure}
	\hspace*{1.9cm}
	\begin{subfigure}[b]{0.35\linewidth}
		\resizebox{1\textwidth}{!}{
			\begin{tikzpicture}
			\begin{axis}[
			xlabel={$l$},ylabel={$\alpha$},
			xmin=0, xmax=1,
			ymin=0, ymax=1.1+2/3,
			axis y line=center, axis x line=center,
			width=10cm, height=10cm,axis on top]
			\addplot[line width=0.25pt,color=white,fill=green!20,opacity=0.75] coordinates{(0,2/3) (0,1+1/3) (2/3,1) (2/3,2/3)} node[rotate=-36,right=-4cm,above=1.05cm,scale=1.3,black,opacity=1] {$\alpha=\frac{5}{3}-l$};
			\addplot[line width=0.25pt,color=white,fill=cyan,opacity=0.25] coordinates{(0,2/3) (0,1+2/3) (1,2/3)} node[rotate=-17,right=-7.0cm,above=-0.3cm,scale=1.3,black,opacity=1] {$\alpha=\frac{4}{3}-\frac{l}{2}$};
			\addplot[line width=1pt,domain=-5:2,samples=2,blue] ({x},{1+2/3-x}); %node[rotate=45,below=0.25cm,left=10cm,black] {$m=\alpha$};
			\addplot[line width=1pt,domain=2/3:1,samples=2,green] ({2/3},{x}); 
			\addplot[line width=1pt,domain=-2:2/3,samples=2,green] ({x},{1+1/3 -x/2}); 
			%%% Linea orizzontale
			%	\addplot[line width=1pt,domain=0:2/3,samples=2,green] ({x},{2/3}); 
			%	\addplot[line width=1pt,domain=2/3:1,samples=2,blue] ({x},{2/3}); 
			%Linea alpha=1
			\addplot[line width=0.1pt,dashed,domain=0:1,samples=2,black] ({x},{1}); 	
			\addplot[line width=0.1pt,dashed,domain=0:1,samples=2,black] ({x},{1+2/3}); 			
			\addplot[line width=0.1pt,dashed,domain=0:2,samples=2,black] ({x},{2/3}); 	
			\end{axis}
			\filldraw[black] (5,6.5) node[anchor=north west] {{\bf{Superlinear}}};
			\filldraw[black] (5,6) node[anchor=north west] {{\bf{chemoattractant}}};			
			\filldraw[black] (5,4+1/2) node[anchor=north west] {{\bf{Sublinear}}};
			\filldraw[black] (5,4) node[anchor=north west] {{\bf{chemoattractant}}};
			%%Parentesi
			%% Parentesi
			\draw
			[decorate,decoration={brace,amplitude=10pt,mirror},xshift=-4pt,yshift=0pt]
			(4.65,4.8) -- (4.65,7.928)node [black,midway,xshift=13pt] {};
			\draw
			[decorate,decoration={brace,amplitude=10pt,mirror},xshift=-4pt,yshift=0pt]
			(4.65,3.2) -- (4.65,4.8)node [black,midway,xshift=13pt] {};
			\end{tikzpicture}
		}
		\caption{Case $n=3$. PP: $l\in (0,\frac{2}{3})$, $\frac{2}{3}\leq \alpha<\frac{4}{3}-\frac{l}{2}$. PE:  $l\in (0,1)$, $\frac{2}{3}\leq \alpha<\frac{5}{3}-l.$} \label{fig:M2}  
	\end{subfigure}
	%%n=4 and n=5 
	\begin{subfigure}[b]{0.35\linewidth}
		\resizebox{1\textwidth}{!}{
			\begin{tikzpicture}
			\begin{axis}[
			xlabel={$l$},ylabel={$\alpha$},
			xmin=0, xmax=1,
			ymin=0, ymax=1.1+2/4,
			axis y line=center, axis x line=center,
			width=10cm, height=10cm,axis on top]
			\addplot[line width=0.25pt,color=white,fill=green!20,opacity=0.75] coordinates{(0,2/4) (0,1+1/4) (2/4,1) (2/4,2/4)} node[rotate=-32,right=-3.3cm,above=2.35cm,scale=1.3,black,opacity=1] {$\alpha=\frac{3}{2}-l$};
			\addplot[line width=0.25pt,color=white,fill=cyan,opacity=0.25] coordinates{(0,2/4) (0,1+2/4) (1,2/4)}
			node[rotate=-24,right=-6.96cm,above=-0.35cm,scale=1.3,black,opacity=1] {$\alpha=\frac{5}{4}-\frac{l}{2}$};
			\addplot[line width=1pt,domain=-5:2,samples=2,blue] ({x},{1+2/4-x}); 
			\addplot[line width=1pt,domain=-2:2/4,samples=2,green] ({x},{1+1/4 -x/2});
			\addplot[line width=1pt,domain=2/4:1,samples=2,green] ({2/4},{x});  
			%%% Linea orizzontale
			%	\addplot[line width=1pt,domain=0:2/4,samples=2,green] ({x},{2/4}); 
			%	\addplot[line width=1pt,domain=2/4:1,samples=2,blue] ({x},{2/4}); 	
			%Linea alpha=1
			\addplot[line width=0.1pt,dashed,domain=0:1,samples=2,black] ({x},{1}); 	
			\addplot[line width=0.1pt,dashed,domain=0:1,samples=2,black] ({x},{1+2/4}); 	
			\addplot[line width=0.1pt,dashed,domain=0:2,samples=2,black] ({x},{2/4}); 					
			\end{axis}
			\filldraw[black] (6.5,7) node[anchor=north west] {{\bf{Superlinear}}};
			\filldraw[black] (6.5,6.5) node[anchor=north west] {{\bf{chemoattractant}}};
			\filldraw[black] (6.5,4+1/2) node[anchor=north west] {{\bf{Sublinear}}};
			\filldraw[black] (6.5,4) node[anchor=north west] {{\bf{chemoattractant}}};
			%Parentesi
			\draw
			[decorate,decoration={brace,amplitude=10pt,mirror},xshift=-4pt,yshift=0pt]
			(6,5.26) -- (6,7.828)node [black,midway,xshift=13pt] {};
			\draw
			[decorate,decoration={brace,amplitude=10pt,mirror},xshift=-4pt,yshift=0pt]
			(6,2.6) -- (6,5.26)node [black,midway,xshift=13pt] {};
			\end{tikzpicture}
		}
		\caption{Case $n=4$. PP: $l\in (0,\frac{1}{2})$, $\frac{1}{2}\leq \alpha<\frac{5}{4}-\frac{l}{2}$. PE:  $l\in (0,1)$, $\frac{1}{2}\leq \alpha<\frac{3}{2}-l$.} \label{fig:}
	\end{subfigure}
	\hspace*{1.9cm}
	\begin{subfigure}[b]{0.35\linewidth}
		\resizebox{1\textwidth}{!}{
			\begin{tikzpicture}
			\begin{axis}[
			xlabel={$l$},ylabel={$\alpha$},
			xmin=0, xmax=1,
			ymin=0, ymax=1.1+2/5,
			axis y line=center, axis x line=center,
			width=10cm, height=10cm,axis on top]
			\addplot[line width=0.25pt,color=white,fill=green!20,opacity=0.75] coordinates{(0,2/5) (0,1+1/5) (2/5,1) (2/5,2/5)} node[rotate=-30,right=-3.3cm,above=3.06cm,scale=1.3,black,opacity=1] {$\alpha=\frac{7}{5}-l$};
			\addplot[line width=0.25pt,color=white,fill=cyan,opacity=0.25] coordinates{(0,2/5) (0,1+2/5) (1,2/5)} node[rotate=-25,right=-7.7cm,above=-0.3cm,scale=1.3,black,opacity=1] {$\alpha=\frac{6}{5}-\frac{l}{2}$};
			\addplot[line width=1pt,domain=-5:2,samples=2,blue] ({x},{1+2/5-x}); 
			\addplot[line width=1pt,domain=-2:2/5,samples=2,green] ({x},{1+1/5 -x/2}); 
			\addplot[line width=1pt,domain=2/5:1,samples=2,green] ({2/5},{x});  
			%%% Linea orizzontale
			%	\addplot[line width=1pt,domain=0:2/5,samples=2,green] ({x},{2/5}); 
			%	\addplot[line width=1pt,domain=2/5:1,samples=2,blue] ({x},{2/5}); 
			%Linea alpha=1
			\addplot[line width=0.1pt,dashed,domain=0:1,samples=2,black] ({x},{1}); 	
			\addplot[line width=0.1pt,dashed,domain=0:1,samples=2,black] ({x},{1+2/5}); 		
			\addplot[line width=0.1pt,dashed,domain=0:2,samples=2,black] ({x},{2/5}); 	
			\end{axis}
			\filldraw[black] (6,7.1) node[anchor=north west] {{\bf{Superlinear}}};
			\filldraw[black] (6,6.6) node[anchor=north west] {{\bf{chemoattractant}}};			
			\filldraw[black] (6,4+1/2) node[anchor=north west] {{\bf{Sublinear}}};	
			\filldraw[black] (6,4) node[anchor=north west] {{\bf{chemoattractant}}};
			%% Parentesi
			\draw
			[decorate,decoration={brace,amplitude=10pt,mirror},xshift=-4pt,yshift=0pt]
			(5.6,5.64) -- (5.66,7.828)node [black,midway,xshift=13pt] {};
			\draw
			[decorate,decoration={brace,amplitude=10pt,mirror},xshift=-4pt,yshift=0pt]
			(5.6,2.2) -- (5.6,5.64)node [black,midway,xshift=13pt] {};
			\end{tikzpicture}
		}
		\caption{Case $n=5$. PP: $l\in (0,\frac{2}{5})$, $\frac{2}{5}\leq \alpha<\frac{6}{5}-\frac{l}{2}$. PE:  $l\in (0,1)$, $\frac{2}{5}\leq \alpha<\frac{7}{5}-l$.} \label{fig:M2}  
	\end{subfigure}
	\caption{Illustration comparing for some values of the dimension $n$ the regions in the $l \alpha$-plane where both parabolic-parabolic (PP, green sector) and parabolic-elliptic (PE, cyan sector) models from problem \eqref{problem} possess uniformly bounded solutions. The superlinear ($\alpha>1$) and sublinear $(\alpha<1)$  chemoattractant zones are also marked. }\label{FigComparisonParParParEll}
\end{figure} 
\subsection{A view to the parabolic-elliptic case}\label{SubSectionOvervieParaEll}
When the parabolic-parabolic problem \eqref{problem} is simplified into parabolic-elliptic, with equation for the chemical replaced by $0=\Delta v-v+g(u)$, assumption \eqref{AssumptionAlphaElleN} becomes sharper; precisely $\frac{2}{n}\leq \alpha <1+\frac{2}{n}-l$, which requires for compatibility only the restriction $l\in (0,1)$. We highlights this aspect in Figure \ref{FigComparisonParParParEll}, where we overlap the regions defined by the interplay between $\alpha$ and $l$ in both models; in particular, we also distinguish the zones with superlinear chemoattractant ($\alpha>1$) and sublinear chemoattractant ($\alpha<1$).  We understand that the observed gap between the range of parameters is not only justified by some technical reasons (see Remark \ref{HintsEllipticCase} at the end of the paper, where few mathematical indications are given) but also by biological ones. Indeed, the fact that in the simplified version  the values of $v(\cdot, t)$ only depend on the values
of $u(\cdot, t)$ at the same time, is a strong modeling assumption.   It corresponds to the situation where the signal responses to the concentration of the particles much faster than the organisms do to the signal; in particular, such difference in the relative adjustment of the bacteria and the chemoattractant makes that this one reaches its equilibrium instantaneously.	
\subsection{Organization of the paper}
The rest of the paper is structured in this way. Section \ref{SectionLocalInTime} is concerned with the local existence question of classical solutions to \eqref{problem} and some of their properties. Some general inequalities are included in Section \ref{MainInequalitySection}. They are mainly devoted to establish how to
ensure globability and boundedness of local solutions using their boundedness in some proper Sobolev spaces; a key cornerstone in this direction is the  procedure to fix the corresponding exponents of theses spaces ($\S$\ref{ProcedureFixParametrSection}).  Finally, the mentioned bound is derived in Section \ref{AprioriEstimatesSection}, which also includes the proof of Theorem \ref{MainTheorem}. 
\section{Existence of local-in-time solutions and main properties}\label{SectionLocalInTime}
Let us dedicate to the existence of classical solutions to system \eqref{problem}. It is show that such solutions are at least local and,  additionally, satisfy some crucial estimates.
\begin{lemma}[\rm{Local existence}]\label{LocalExistenceLemma}   
	%{\textcolor{blue}{
	Let $\Omega$ be a bounded and smooth domain of $ \R^n$, with $n \geq 2$.  Moreover, let $f$ and $g$ fulfill \eqref{f} and \eqref{AssumptionOnSegregationG}, respectively, with $l\in (0,\frac{2}{n})$ and $\alpha$ satisfying \eqref{AssumptionAlphaElleN}. Then, for any  nontrivial $(u_0,v_0)\in C^0(\bar{\Omega})\times C^{1}(\bar\Omega)$, with $u_0\geq 0$ and $v_0\geq 0$ on $\bar{\Omega}$, there exists $\TM \in (0,\infty]$ and a unique pair of nonnegative functions $(u,v)\in 
 (C^0(\bar{\Omega}\times [0,\infty))\cap  C^{2,1}(\bar{\Omega}\times (0,\infty)))^2,
	$
%	$$
%\begin{aligned}
%u &\in  C^0(\bar{\Omega}\times [0,\TM))\cap  C^{2,1}(\bar{\Omega}\times (0,\TM)),\\
%v &\in  C^0(\bar{\Omega}\times [0,\TM))\cap  C^{2,1}(\bar{\Omega}\times (0,\TM))\cap L_{loc}^\infty([0,\TM);W^{1,r}(\Omega)),
%\end{aligned}
%$$
	such that this dichotomy criterion holds true:
\begin{equation}\label{ExtensibilityCrit}
 \textrm{either}\,\; T_{max}=\infty\; \textrm{or}\; \limsup_{t \rightarrow T_{max}}(\lVert u(\cdot,t)\rVert_{L^\infty(\Omega)}+\lVert v(\cdot,t)\rVert_{L^{\infty}(\Omega)}) =\infty. 
\end{equation} 
In addition, the $u$-component obeys the mass conservation property, i.e. 
\begin{equation}\label{massConservation}
	\int_\Omega u(x, t)dx =\int_\Omega u_0(x)dx=m>0\quad \textrm{for all }\, t \in (0,\TM),
\end{equation}
whilst for some $c_0>0$ the $v$-component is such that 
\begin{equation}\label{Cg}
\lVert v (\cdot, t)\rVert_{W^{1,n}(\Omega)}\leq c_0 \quad \textrm{on } \,  (0,\TM).
\end{equation}
\begin{proof}
We just mention that the conclusions concerning the  local-in-time well-posedness as well as the dichotomy criterion \eqref{ExtensibilityCrit} can be established by straightforward adaptations of widely used 
methods involving an appropriate fixed point framework and standard parabolic regularity theory; in particular we cite \cite[Theorem 3.1]{HorstWink} for the case $g(u)=u$ and \cite[Lemma 3.1]{ViglialoroWoolleyAplAnal} for $g$ as in our hypotheses. Moreover, comparison arguments apply to yield both $u,v\geq 0$ in $\Omega\times (0,\TM)$. 

On the other hand, the mass conservation property easily comes by integrating over $\Omega$ the first equation of \eqref{problem}, in conjunction with the boundary and initial conditions.

Finally, the last claim is derived as follows. 	From the assumption $0<l<\frac{2}{n}$, we can first of all fix $\frac{n}{2}<\gamma<n$ complying with $\gamma\leq \frac{1}{l}$. In this way, thorough the Hölder inequality, taking in mind
\eqref{AssumptionOnSegregationG} and the mass conservation property \eqref{massConservation}, we have
\begin{equation}\label{Estim_general_For_u^p}
\int_\Omega g(u)^\gamma\leq  K_0^\gamma \int_\Omega u^{\gamma l}\leq K_0^\gamma  m^{\gamma l}|\Omega|^{1-\gamma l}\quad \textrm{for all } t< \TM.
\end{equation}		
Henceforth, we can also  pick $\frac{1}{2}<\rho <1$ such that $\zeta=1-\rho -\frac{n}{2}(\frac{1}{\gamma}-\frac{1}{n})>0$. Subsequently, since by means of the representation formula for $v$ we have
\[v(\cdot,t) =e^{t(\Delta-1)}v_0 +\int_0^t e^{(t-s)(\Delta-1)}g(u(\cdot,s))ds \quad \textrm{for all } \,t\in (0,T_{max}),
\]
we invoke properties regarding the Neumann heat semigroup $(e^{t \Delta})_{t \geq 0}$ (see Section 2 of \cite{HorstWink} and Lemma 1.3 of \cite{WinklAggre}), so to obtain for some $\lambda_1>0$ and $C_S>0$
\begin{equation}\label{Bound_v_1-qA}
\begin{split}
 \lVert  v (\cdot, t) \lVert_{W^{1,n}(\Omega)} & \leq  e^{-t}  \lVert e^{t \Delta}v_0 \lVert_{W^{1,n}(\Omega)}+\int_0^t \lVert  e^{(t-s)(\Delta-1)} g(u(\cdot,s))\lVert_{W^{1,n}(\Omega)}ds \\ &
\leq C_S \lVert v_0\lVert_{W^{1,n}(\Omega)}+C_S\int_0^t \lVert  (-\Delta +1)^\rho e^{(t-s)(\Delta-1)} g(u(\cdot, s))\lVert_{L^n(\Omega)}ds
\\ &
\leq C_S \lVert v_0\lVert_{W^{1,n}(\Omega)}+C_S\int_0^t (t-s)^{-\rho  -\frac{n}{2}(\frac{1}{\gamma}-\frac{1}{n})}e^{-\lambda_1 (t-s)}\lVert g(u(\cdot, s))\lVert_{L^\gamma(\Omega)}ds. 
%\\ &
%\leq d_1 \lVert \nabla v_0\lVert_{L^\infty(\Omega)}+d_2 \Big(\int_0^t (1+(t-s)^{-\frac{1}{2}-\frac{n}{2p}}ds\Big)^\frac{p}{p-1} \Big(\int_0^t\lVert q(\cdot, %s)\lVert_{L^p(\Omega)}^p ds\Big)^\frac{1}{p}  \\ & 
%\leq d_1 \lVert \nabla v_0\lVert_{L^\infty(\Omega)}+d_2 d_4 C occhio!!!!!!!.
%	\leq C_S \lVert v_0\lVert_{W^{1,n}(\Omega)}+C_SC(p)\int_0^t (t-s)^{-\rho  -\frac{n}{2}(\frac{1}{\gamma}-\frac{1}{n})}e^{-\lambda_1 (t-s)}ds, 
\end{split}
\end{equation}
As a consequence, the introduction of the Gamma function $\Gamma$ infers
\begin{equation*}
\int_0^t (t-s)^{-\rho  -\frac{n}{2}(\frac{1}{\gamma}-\frac{1}{n})}e^{-\lambda_1 (t-s)}ds\leq \lambda_1^{-\zeta} \Gamma(\zeta),
\end{equation*}
which combined with  bounds \eqref{Estim_general_For_u^p} and \eqref{Bound_v_1-qA}  conclude the proof.  
\end{proof}
\end{lemma}
In the sequel of the paper with $(u,v)$ we will refer to the gained local classical solution to problem \eqref{problem}; in particular we might tacitly avoid to mention that such solution is produced by the initial data  $(u_0,v_0)$.
\section{Preliminaries: inequalities and parameters}\label{MainInequalitySection}
With the local solution $(u,v)$ to  problem \eqref{problem} at disposal,  its uniform boundedness on $(0,\TM)$ is achieved when uniform-in-time bound for $u$ in some $L^p$-norm and for $|\nabla v|^2$ in some $L^q$-norm, with proper $p$ and $q$, is derived. This will be obtained by constructing an absorption inequality satisfied by the functional 
\begin{equation}\label{FunctionalDefinedOnLocalSol}
y(t):=\frac{1}{p(p-1)}\int_\Omega (u+1)^p+\frac{1}{q}\int_\Omega |\nabla v|^{2q} \quad \textrm{for all } t\in (0,\TM).
\end{equation}
 In particular, the entire procedures requires, first, to adequately  manipulate inequalities resulting by differentiating $y(t)$ with respect to the time and, secondly, to figure out how to choose the parameters $p$ and $q$. In view of its decisive role, this second step will be discussed with some details after this subsection.
\subsection{Some algebraic and functional inequalities}
This three coming lemmas will be used in the next logical steps. We start by considering a suitable version of the Gagliardo--Nirenberg interpolation inequality, commonly used to treat nonlinearities appearing
in the diffusion and/or chemosensitivity terms (see \cite{TaoWangEtALDCDS-B,ViglialoroAppliedMathOpt,WangYangYuDCDS-A}), and successively by recalling a particular boundary integral used to deal with terms defined in non-convex domains.
\begin{lemma}[\rm{Gagliardo--Nirenberg inequality}]\label{LemmaGN}
	Let $\Omega$ be a bounded and smooth domain of $\mathbb{R}^n$, with $n\geq 1$, and  $0<\mathfrak{q}\leq \mathfrak{p} \leq \infty$  satisfying	
	$
	\frac{1}{2} \leq \frac{1}{n} + \frac{1}{\mathfrak{p}}. 
	$
	Then, for $
	a= \frac{\frac{1}{\mathfrak{q}}-\frac{1}{\mathfrak{p}}}{\frac{1}{\mathfrak{q}}+\frac{1}{n}-\frac{1}{2}},
	$
	there exists $C_{GN}=C_{GN}(\mathfrak{p},\mathfrak{q},\Omega)>0$ such that
	\[
	\|w\|_{L^{\mathfrak{p}}(\Omega)} \leq C_{GN} (\|\nabla w\|_{L^2(\Omega)}^a \|w\|_{L^{\mathfrak{q}}(\Omega)}^{1-a}+ \|w\|_{L^{\mathfrak{q}}(\Omega)}) 
	\quad \text{ for all } w \in W^{1,2}(\Omega) \cap L^{\mathfrak{q}}(\Omega).
	\]
	\begin{proof}
		See \cite[Lemma 2.3]{LiLankeitNonLinearity}.
	\end{proof}
\end{lemma}
\begin{lemma}\label{LemmaConv}
	Let $\Omega$ be a bounded and smooth domain of  $\mathbb{R}^n$, with $n\geq 1$, and $q \in [1, \infty)$. Then for any $\eta >0$ there is $C_{\eta}>0$ such that for any 
	$w \in C^2(\overline{\Omega})$ satisfying $\frac{\partial w}{\partial \nu} = 0$ on $\partial \Omega$ this inequality holds:
	\[
	\int_{\partial \Omega} |\nabla w|^{2q-2} \frac{\partial |\nabla w|^2}{\partial \nu} \leq \int_{\Omega} |\nabla|\nabla w|^q|^2 + C_{\eta}.
	\]
	\begin{proof}
		A proof can be found in \cite[Propostion 3.2]{YokotaEtAlNonCONVEX}. It is the combination of a proper Sobolev embedding, an estimate of the boundary term as an expression depending on the curvature of $\partial \Omega$ and $|\nabla w|^2$ as well as  applications of a fractional Gagliardo--Nirenberg inequality  and Young's one. (When $\Omega$  is convex, the left hand side of the inequality is nonpositive (see \cite[Appendix]{DalPassoEtAlSIAM} and \cite[Lemma 3.2]{TaoWinkParaPara}.)
	\end{proof}
\end{lemma}
Thanks to the next results (variants of Young’s inequality), products of powers will be estimated by suitable sums involving their bases and powers of sums controlled by sums of powers. 
\begin{lemma} \label{LemmaIneq}
	Let $a,b \geq 0$ and $d_1, d_2>0$ such that $d_1 + d_2 <1$. Then for all $\epsilon >0$ there exists $c>0$ such that
	\[
	a^{d_1}b^{d_2} \leq \epsilon(a+b) +c. 
	\]
	Moreover, for further  $d_3, d_4>0$, it is possible to find positive $d_5$ and $d$ such that
	\[
	a^{d_3} + b^{d_4} \geq 2^{-d_5} (a+b)^{d_5} -d.
	\]
	\begin{proof}
		We show the first inequality, being the proof of the second similar. By applying Young's inequality with conjugate exponents $\frac{1}{d_1}$ and $\frac{1}{1-d_1}$, we obtain for any $\epsilon_1>0$ and some $c_1(\epsilon_1)>0$ that
		\[
		a^{d_1}b^{d_2}= a^{d_1}(b^{\frac{d_2}{1-d_1}})^{1-d_1} \leq \epsilon_1 a+ c_1(\epsilon_1) b^{\frac{d_2}{1-d_1}}.
		\]
		Moreover, due to $\frac{d_2}{1-d_1} <1$, a further application of the same inequality to the latter power provides for every positive $\epsilon_2$ and proper $c_2(\epsilon_2)>0$ this relation: $
		c_1(\epsilon_1) b^{\frac{d_2}{1-d_1}} \leq \epsilon_2 b + c_2(\epsilon_2).$	By putting together the two inequalities and choosing $\epsilon_1 = \epsilon_2$, the first part of the lemma is concluded. All the details of the second inequality can be found in \cite[Lemma 3.3]{MarrasViglialoroMathNach}.
	\end{proof}
\end{lemma}
\subsection{The right procedure in fixing the parameters $p$ and $q$}\label{ProcedureFixParametrSection}
In this sequence of lemmas, we will verify that the mutual relation between the parameters $l$ in \eqref{f} and $\alpha$ and $l$ in \eqref{AssumptionOnSegregationG}, i.e. relation \eqref{AssumptionAlphaElleN}, is such that the mentioned parameters $p$ and $q$ may be chosen in the appropriate way to make sure our general machinery work. 
\begin{lemma}\label{LemmaPreparingRelationAphaElleandN}
	For any $n \in \N$, with $n\geq 2$, let $l\in (0,\frac{2}{n})$ and $\alpha$ comply with assumption \eqref{AssumptionAlphaElleN}.
	Then there exist  $1<\theta<\frac{n}{n-2}$ and  $\mu>\frac{n}{2}$ such that 
	\begin{equation}\label{RelationAlphaEllen>2}
	\frac{ l (2\mu-1)}{4\mu-n}<\frac{n(\theta +1-2\alpha \theta)+2\theta}{2n\theta +n^2-n^2\theta}.
	\end{equation}
%	
%	
%	following conclusions hold true:
%	\begin{enumerate}[label=\roman*.)]
%		\item\label{FirtItemRelationAlphan=2} for $n=2$, there exists $\theta>1$ such that 
%		\begin{equation}\label{RelationAlphaEllen=2}
%		\frac{1}{2\theta-2\alpha\theta+1}<\frac{1}{l},
%		\end{equation}
%		\item \label{SecondItemRelationAlphan>2} for $n\geq 3$, there exist $1<\theta<\frac{n}{n-2}$ and  $\mu>\frac{n}{2}$ such that 
%			\begin{equation}\label{RelationAlphaEllen>2}
%	\frac{ l (2\mu-1)}{4\mu-n}<\frac{n(\theta +1-2\alpha \theta)+2\theta}{2n\theta +n^2-n^2\theta},
%	\end{equation}
%%		\item \label{ThirdItemRelationAlpha=1n>2} for $n\geq 3$, there exists $\mu>\frac{n}{2}$ such that for  $\alpha=1$ 
%%	\begin{equation}\label{RelationAlphaEllen>2ForAlpha=1}
%%	\frac{ l (2\mu-1)}{4\mu-n}<\frac{1}{n}.
%%		\end{equation}
%	\end{enumerate}
\begin{proof}
	From now on, we first precise that for $n=2$, $1<\theta <\frac{n}{n-2}$ indicates that $\theta$ might also be fixed large as we want; despite that, this is not the case, and we will take $\theta$ always sufficiently close to 1. 
	
Precisely, for $n\geq 2$,  $1<\theta<\frac{n}{n-2}$ implies that $2n\theta +n^2-n^2\theta>0$; in particular,  we can consider $\theta>1$ small enough so to have $n(\theta +1-2\alpha \theta)+2\theta>0$. Hence the function 
$$h(\theta,\mu):=	\frac{ l (2\mu-1)}{4\mu-n}-\frac{n(\theta +1-2\alpha \theta)+2\theta}{2n\theta +n^2-n^2\theta},$$
is the difference of two positive terms. Now, from our assumptions
$$\lim_{\mu\to +\infty} h(1,\mu)= \Big(\frac{l}{2}-1+\alpha-\frac{1}{n}\Big)<0,$$
and the claim  is proved by means of continuity arguments. 
%Finally, the claim in item \ref{ThirdItemRelationAlpha=1n>2} follows from the observation that $h(\theta,\mu; \alpha)$ is for $\alpha=1$ a $\theta$-independent function. More precisely $h(\theta,\mu; 1)=	\frac{ l (2\mu-1)}{4\mu-n}-\frac{1}{n}$ converges to $\frac{l}{2}-\frac{1}{2}$ which is negative in view of $l<\frac{2}{n}.$ So, from sign-preserving properties for continuous functions, we conclude that the same definitively holds for some $\mu>\frac{n}{2}$ sufficiently large. 

\end{proof}
\end{lemma}

\begin{lemma}\label{pqtildeLemman=2andn>2Lemma}  
	Let the hypotheses of Lemma \ref{LemmaPreparingRelationAphaElleandN} be satisfied, and $1<\theta <\frac{n}{n-2}$ and $\mu>\frac{n}{2}$ be therein fixed. Then there is $q_r\in [1,\infty)$ such that for all $q>q_r$ one has this compatibility relation:
	\begin{equation}\label{ptildeqtilden>2}
	-\frac{2l (nq+\mu(4+n^2-2n(2+q)))}{n(4\mu-n)}= f_1(q)<f_2(q)=\frac{q(2n(\theta +1-2\alpha \theta)+4\theta)+2n\theta (\alpha-1)(n-2)}{2n\theta +n^2-n^2 \theta}.
	\end{equation}
%whereas for $n\geq 3$, $\alpha=1$, and some $\mu>\frac{n}{2}$,
%\begin{equation}\label{ptildeqtilden>2ForAlpha=1}
%-\frac{2l (n\tilde{q}+\mu(4+n^2-2n(2+\tilde{q})))}{n(4\mu-n)}	<\tilde{p}<\frac{2\tilde{q}}{n}.
%\end{equation}
\begin{proof}
 Some easy computations show that for any $n\geq 2$ the claim follows once it is established that for $1<\theta<\frac{n}{n-2}$ and $\mu>\frac{n}{2}$ as in the hypotheses, and for 
\begin{align}
\mathcal{A}  = \mathcal{A}(\theta,\mu)&= 	\frac{2 l (2\mu-1)}{4\mu-n},            &  \mathcal{C}  = \mathcal{C}(\theta,\mu) &=\frac{2n(\theta +1-2\alpha \theta)+4\theta}{2n\theta +n^2-n^2\theta},  \nonumber\\ 
\mathcal{B}  = \mathcal{B}(\theta,\mu) &= \frac{2l\mu(n-2)^2}{n(4\mu-n)},  &  \mathcal{D}  = \mathcal{D}(\theta,\mu)  &= \frac{2n\theta (1-\alpha)(n-2)}{2n\theta +n^2-n^2\theta},\nonumber
\end{align}
there is a  $q\geq 1$ entailing
\begin{equation}\label{RelationonCapitalABCDForqTilde}
	\mathcal{A}-\frac{\mathcal{B}}{q}<	\mathcal{C}-\frac{\mathcal{D}}{q}.
\end{equation}  
As we  justify and explain it in Figure \ref{FigureForQTilde}, the above occurs whenever  $\mathcal{A}-\mathcal{C}<0$.
\end{proof}
\end{lemma}
 \begin{SCfigure}[1.3][htb]
	\centering
	\begin{tikzpicture}
	\begin{axis}[
	axis lines=middle,
	samples=200,
	xtick={1,...,8},
	ytick={-2,...,10}
	]
	\addplot[cyan,domain=0.6:8.5] {3/(x) -2};
	\addplot[line width=0.25mm,green,dashed,domain=0.15:8.5] {-2};
	\draw (140,169)  node  {${q}_{r}$};
	\draw (796,219)  node  {$q$};
	\draw (22,219)  node  {$O$};
	\draw (95,475)  node  {$k(\theta, \mu; q)$};
	\draw[fill=blue] (135,200) circle(0.08211cm);
	\draw [thick,dashed] (133,218)-- (133,255);
	\draw [thick,dashed] (133,255)-- (438,255);
	\draw [thick,dashed] (525,255)-- (800,255);
	\draw (480,255)  node  {$q$};
	\draw[fill=white] (135,255) circle(0.08211cm);
	%\draw[red!20,dashed] (axis cs:0,-1.3) -- (axis cs:4,-1.3);
	\end{axis}
	\end{tikzpicture}
	\caption{By setting $k(\theta,\mu; q):=\mathcal{A}-\mathcal{C}-\frac{\mathcal{B}}{q}	+\frac{\mathcal{D}}{q}$, from inequality \eqref{RelationonCapitalABCDForqTilde} we intend to find $q_r$ such that for some $\theta$ and $\mu$ we have that  $k(\theta,\mu; q)<0$ for all $q>q_r$. Let $\theta$ close to $1$ and  $\mu$ sufficiently large be taken from  Lemma \ref{LemmaPreparingRelationAphaElleandN}. As a consequence, inequality \eqref{RelationAlphaEllen>2} leads to $\mathcal{A}-\mathcal{C}<0$, whereas $\mathcal{B}-\mathcal{D}\in \R$. In particular, considering that $\frac{\partial k(\cdot, \cdot; q)}{\partial q}=\frac{\mathcal{B}-\mathcal{D}}{q^2}$
		and $
		\lim_{q\to +\infty} k(\cdot, \cdot; q)= \mathcal{A}-\mathcal{C},
		$
		the illustration shows the qualitative behavior of the function $k(\theta,\mu; q)$ for these values of $\theta$ and $\mu$, assuming the nontrivial situation $\mathcal{B}-\mathcal{D}<0$. (If, indeed, $\theta$ and $\mu$ infer $\mathcal{B}-\mathcal{D}\geq 0$, $k$ is negative for all $q$.) Then, by indicating  with $q_{r}=$$\frac{\mathcal{B}-\mathcal{D}}{\mathcal{A}-\mathcal{C}}$ the root of $k$, any $q\in (q_{r},\infty)$ satisfies relation \eqref{RelationonCapitalABCDForqTilde}. (In order to clarify the role of $\mu$ and $\theta$, we observe that for $n=2$ the chain of inequality in \eqref{ptildeqtilden>2} is more manageable; in fact, it reads $
		l q<q(2\theta-2\alpha\theta+1),
		$ directly coming from
		$
		l<	2\theta-2\alpha\theta+1,
		$ corresponding to \eqref{RelationAlphaEllen>2} when $n=2$, and it is $\mu$-independent and true for some $\theta$ approaching 1, once $\alpha<\frac{3-l}{2}$ from \eqref{AssumptionAlphaElleN} is considered.)}	\label{FigureForQTilde}
	%summarizing properties of classical solutions to problem \eqref{problem} when arbitrary $\chi>0$ and $\lambda \in \R$ are considered. The blue line stands for global and uniformly bounded solutions (\cite[Corollary 2.6]{TelloWinkParEl}); the green line, for global solutions (Theorem \ref{MainTheorem}) and the red line, for blow-up solutions ()\cite[Th. 1.1]{WinklerBlowUpLowDimensionWithLogistic}). No result is available for values of $k$ belonging to the gray intervals. (The result in the horizontal cyan line comes is a consequences of standard derivations.)}
\end{SCfigure}
\begin{lemma}\label{LemmaCoefficientAiAndExponents} 
	Let the hypotheses of Lemma \ref{LemmaPreparingRelationAphaElleandN} be satisfied,  and $1<\theta <\frac{n}{n-2}$ and $\mu>\frac{n}{2}$ be therein fixed. Then there are $p\in[1,\infty)$ and $q\in[1,\infty)$ such that 
%	\begin{equation}\label{CoefficienAeC_Gagl-Nir}
%	a_1=\frac{\frac{np}{2}(1-\frac{1}{(p+2\alpha-2)\theta})}{1-\frac{n}{2}+\frac{np}{2}} \quad \quad \quad  a_2=\frac{nq(\frac{1}{n}-\frac{1}{2\theta'})}{1-\frac{n}{2}+q}
%	\end{equation}
%		\begin{equation}\label{CoefficienBeD_Gagl-Nir}
%	a_3=\frac{\frac{np}{2l}(1-\frac{1}{2\mu})}{1-\frac{n}{2}+\frac{np}{2l}} \quad \quad \quad  a_4=\frac{nq(\frac{1}{n}-\frac{1}{2(q-1)\mu'})}{1-\frac{n}{2}+q}
%	\end{equation}
	\begin{align}
	a_1&= 	\frac{\frac{np}{2}(1-\frac{1}{(p+2\alpha-2)\theta})}{1-\frac{n}{2}+\frac{np}{2}},            &  a_2&=\frac{nq(\frac{1}{n}-\frac{1}{2\theta'})}{1-\frac{n}{2}+q},  \nonumber \\ 
a_3 &= \frac{\frac{np}{2l}(1-\frac{1}{2\mu})}{1-\frac{n}{2}+\frac{np}{2l}},  &  a_4  &= \frac{nq(\frac{1}{n}-\frac{1}{2(q-1)\mu'})}{1-\frac{n}{2}+q},\nonumber\\ 
\kappa_1 &  =\frac{\frac{np}{2} (1- \frac{1}{p})}{1-\frac{n}{2}+\frac{np}{2}} \nonumber, & \kappa_2 & =  \frac{q - \frac{n}{2}}{1-\frac{n}{2}+q}, 
	\end{align}
belong to the interval $(0,1)$ and, additionally, imply that these other relations hold true: 
		\begin{equation}\label{MainInequalityExponents}
		%\begin{cases}
	\beta_1 + \gamma_1 =\frac{p-2+2\alpha}{p}a_1+\frac{1}{q}a_2\in (0,1) \;\textrm{ and }\;	\beta_2 + \gamma_2= \frac{2 l }{p}a_3+\frac{q-1}{q}a_4 
	\in (0,1).
		%\end{cases} 
%		\quad \textrm{if }\;\,	\begin{cases}
%		n=2 \, \textrm{ and for some } \theta>1 \textrm{ and any } \mu>1,\\  
%		n\geq 3 \, \textrm{ and for some } \theta<\frac{n}{n-2} \textrm{ and  } \mu>\frac{n}{2}.\\  
%		\end{cases}
		\end{equation}
	\begin{proof}
		For $\theta$ and $\mu$ as in our hypotheses, their conjugate exponents $\theta'$ and $\mu'$ satisfy $\theta' > \frac{n}{2}$ and $
\mu'< \frac{n}{n-2}$. Now, for $q_r$ taken from Lemma \ref{pqtildeLemman=2andn>2Lemma},  since  $(f_1(q),f_2(q))$ is not empty thanks to  compatibility \eqref{ptildeqtilden>2}, we can always consider $p$ and $q$ fulfilling
\begin{equation}\label{Prt_q}
	\begin{cases}
	 q > \max \left\{\frac{(n-2)}{n}\theta', \frac{n}{2\mu'}+1, \frac{2n\theta (\alpha-1)(2-n)}{2n(\theta+1-2\alpha\theta)+4\theta}, q_r\right\} \\
	p>\max \left\{2 + \frac{1}{\theta},\frac{2(n-2)l\mu}{n},\frac{2\theta(\alpha-1)(n-2)}{n-\theta(n-2)}\right\} 
	\end{cases}
\quad \textrm{ and also complying with  } p\in(f_1(q),f_2(q)).
	\end{equation}
Our aim is to show that such restrictions suffice to prove the claim. Straightforward reasoning justify that some of the first relations in  \eqref{Prt_q} imply  $a_1, a_2, a_3, a_4, \kappa_1, \kappa_2 \in (0,1).$ The remaining two inequalities in \eqref{MainInequalityExponents} are, conversely, less direct. Indeed, if it can be immediately inferred that $\frac{p-2+2\alpha}{p}a_1+\frac{1}{q}a_2$ and $\frac{2 l }{p}a_3+\frac{q-1}{q}a_4$ are positive, the other bound requires tedious computations involving $f_1(q)$ and $f_2(q)$. More exactly, algebraic rearrangements give 
\begin{equation}\label{PassoUnoCalcoliPrimaSommaEsponenti<1}
\frac{p-2+2\alpha}{p}a_1+\frac{1}{q}a_2-1=\frac{n^2 (2 (\alpha-1) \theta+p (\theta-1))+2 n (q (-2 \alpha \theta+\theta+1)-\theta(2 \alpha+p-2))+4 q \theta}{\theta (n (p-1)+2) (n-2 (q+1))},
\end{equation}
and 
\begin{equation}\label{PassoUnoCalcoliSecondaSommaEsponenti<1}
\frac{2 l }{p}a_3+\frac{q-1}{q}a_4 -1=\frac{n p (n-4 \mu)-2 l \left(\mu \left(n^2-2 n (q+2)+4\right)+n q\right)}{\mu (n-2 (q+1)) (l(n-2)-n p)}.
\end{equation}
To see that expression \eqref{PassoUnoCalcoliPrimaSommaEsponenti<1} is negative, we notice from the constrains on $p,q,\theta$ and $\mu$ that the denominator is negative, so by  imposing
$$n^2 (2 (\alpha-1) \theta+p (\theta-1))+2 n (q (-2 \alpha \theta+\theta+1)-\theta(2 \alpha+p-2))+4 q \theta>0,$$ we obtain
\begin{equation}\label{A0A0}
p(n^2\theta-n^2-2n\theta)>4\alpha \theta n-4n\theta-4q\theta-2n^2\theta (\alpha-1)-2nq (\theta+1-2\alpha\theta).
\end{equation}
This, taking into account the negativity of $n^2\theta-n^2-2n\theta$, is equivalent to find $q$ such that 
$$4\alpha \theta n-4n\theta-4q\theta-2n^2\theta (\alpha-1)-2nq (\theta+1-2\alpha\theta)<0 \quad \textrm{ or also }\quad 
q>\frac{2n\theta (\alpha-1)(2-n)}{2n(\theta+1-2\alpha \theta)+4\theta},
$$
which is fulfilled by virtue of the choice on $q$ and since the considered  $\theta$ complies with $2n(\theta+1-2\alpha)+4\theta>0$. Subsequently, from \eqref{A0A0} we have that 
\begin{equation}\label{FirstInterplayP-Q}
p<\frac{q(2n(\theta+1-2\alpha\theta)+4\theta)+2n\theta (\alpha-1)(n-2)}{n^2+2n\theta-n^2\theta}
\end{equation}
is satisfied for $p$ and $q$ as in \eqref{Prt_q}.

Let us now turn our attention to \eqref{PassoUnoCalcoliSecondaSommaEsponenti<1}. Unlike the previous case, we immediately see that the denominator is positive and, again by invoking \eqref{Prt_q}, it holds that   
\begin{equation}\label{SecondInterplayP-Q}
p>-\frac{2l (nq+\mu(4+n^2-2n(2+q)))}{n (4\mu-n)}.
\end{equation}
%Since by construction $f_1(q)<p<f_2(q)$, we have compatibility between relations \eqref{FirstInterplayP-Q} and \eqref{SecondInterplayP-Q}; this concludes the proof. 
\end{proof}
\end{lemma}
\begin{remark}\label{RemarkEnlargingPandQ} 
Let us spend some words on how to treat the introduced parameters $p$ and $q$ in accordance with our overall purposes. This technical detail makes the analysis of the present work different and in some sense more thorough with respect those presented in many references above mentioned; therein, indeed, no undesired smallness assumption on $p$, generally, appears. 
\begin{enumerate}[label=(\roman*)]
\item \label{ItemUponEnlarging}Taking the ``lower extremes'' for $q$ in $(q_r,\infty)$ and for $p$ in $(f_1(q),f_2(q))$,  as specified in Lemma \ref{LemmaCoefficientAiAndExponents}, might not be appropriate when dealing with other computations where they are involved. In particular, as we will perform in the last step toward the proof of Theorem \ref{MainTheorem}, it could be necessary to enlarge each one of this values  in order to ensure the validity of certain inequalities/inclusions. Despite that, we understand that some care is needed when this procedure has to be adopted; indeed,  $p$ cannot be taken large as we want independently by $q$, but  this is possible when the order $f_1(q)<p<f_2(q)$ related to relation \eqref{ptildeqtilden>2} is preserved. (This was already imposed in the same Lemma  \ref{LemmaCoefficientAiAndExponents}.) 
\item In support to the previous item,  we point out that even though asymptotically we have
\begin{equation*}
\begin{cases}
\frac{p-2+2\alpha}{p}a_1=\frac{n (\theta (2 \alpha+p-2)-1)}{\theta (n (p-1)+2)}\nearrow 1 & \textrm{increasing with } \;p,\\
\frac{1}{q}a_2=\frac{n-2 \theta'}{\theta' (n-2 (q+1))}\nearrow 0& \textrm{decreasing  with } \;q,
\end{cases}
%\end{equation}
\quad and \quad  
%\begin{equation*}
\begin{cases}
\frac{2 l }{p}a_3=\frac{l (1-2 \mu) n}{\mu (l (n-2)-n p)} \nearrow 0 & \textrm{decreasing  with } \;p,\\
\frac{q-1}{q}a_4=\frac{n-2 \mu' (q-1)}{\mu' (n-2 (q+1))}\nearrow 1 & \textrm{increasing  with } \;q,
\end{cases}
\end{equation*}
this is not sufficient to ensure that there exists a couple $(p,q)$ for which both $\frac{p-2+2\alpha}{p}a_1+
\frac{1}{q}a_2<1$ and $\frac{2 l }{p}a_3+
\frac{q-1}{q}a_4<1$ are satisfied. Surely each one of this inequality holds true for two different couples, let's say $(p_0,q_0)$ and $(p_1,q_1)$, but the identification of a single $(p,q)$ producing simultaneously those inequalities requires the extra condition $p\in(f_1(q),f_2(q))$, intimately linked to the main assumption \eqref{AssumptionAlphaElleN}.
\end{enumerate}
\end{remark}
\section{Deriving uniform-in-time $L^p\times L^q$-bounds for $(u,|\nabla v|^2)$. Proof of the main result}\label{AprioriEstimatesSection}
The coming lemma provides a uniform-in-time bound on $(0,\TM)$ for $u$ in $L^p(\Omega)$ and for $|\nabla v|^2$ in $L^q(\Omega)$.
\begin{lemma}\label{LemmaStime}
	Under the hypotheses of Lemma \ref{LocalExistenceLemma} we have the following conclusion: For some $p \in(1,\infty)$ and $q \in (1, \infty)$ there exists $L>0$ such that
	\begin{equation*}% \label{St.u}
	\lVert u(\cdot, t)\|_{L^p(\Omega)} +\|\nabla v(\cdot,t)\|_{L^{2q}(\Omega)}  \leq L \quad \text{ for all } \,t <T_{max}.
	\end{equation*}
	\begin{proof}
		With $\theta$, $\mu$, $p$ and $q$ as in Lemma \ref{LemmaCoefficientAiAndExponents}, the validity of all the computations along this lemma is justified. 
		
		As announced, let us differentiate with respect to the time $y(t)$ defined in  \eqref{FunctionalDefinedOnLocalSol} and split the resulting derivations in three main steps, altogether yielding the proof.
		\subsubsection*{{\bf{Estimating $\frac{1}{p(p-1)}\frac{d}{dt} \int_{\Omega} (u+1)^p$ on $(0,T_{max})$.}}}
		We take $\frac{(u+1)}{p-1}^{p-1}$ as test function for the first equation in \eqref{problem}, so that by integrating by parts we obtain, also in view of the no-flux boundary conditions, that
		\begin{align}\label{Der.u}
		\frac{1}{p(p-1)}\frac{d}{dt} \int_{\Omega} (u+1)^p &= \frac{1}{p-1}\int_{\Omega} (u+1)^{p-1} \nabla \cdot \nabla u - \frac{1}{p-1}
		\int_{\Omega} (u+1)^{p-1} \nabla \cdot (f(u) \nabla v) \\ \nonumber
		%&= - \int_{\Omega} \phi''(u) |\nabla u|^2 + \int_{\Omega} \phi''(u) f(u) \nabla u \cdot \nabla v \\ \nonumber
		&= - \int_{\Omega} (u+1)^{p-2} |\nabla u|^2 + \int_{\Omega} (u+1)^{p-2} f(u) \nabla u \cdot \nabla v \quad \textrm{ on }\, (0, T_{max}).
		\end{align}
		Through  an application  of Young's inequality and \eqref{f},  the latter term reads
		\begin{equation} \label{Y1}
		\int_{\Omega} (u+1)^{p-2} f(u) \nabla u \cdot \nabla v \leq \frac{1}{2} \int_{\Omega} (u+1)^{p-2} |\nabla u|^2 + \frac{K^2}{2} \int_{\Omega} (u+1)^{p+2\alpha-2} |\nabla v|^2 \quad \textrm{ for all } \,t \in (0, T_{max}),
		\end{equation}
		and  the second integral at the right-hand side is estimated by the H\"{o}lder inequality so to have 
		\begin{equation} \label{H1}
		\frac{K^2}{2}\int_{\Omega} (u+1)^{p+2\alpha-2} |\nabla v|^2 \leq \frac{K^2}{2} \left(\int_{\Omega} (u+1)^{(p+2\alpha-2)\theta}\right)^{\frac{1}{\theta}} 
		\left(\int_{\Omega} |\nabla v|^{2 \theta'}\right)^{\frac{1}{\theta'}} \quad \textrm{ on } (0, T_{max}).
		\end{equation}
		Now we can apply Lemma \ref{LemmaGN} with $\mathfrak{p}=\frac{2(p-2+2\alpha)\theta}{p}$, $\mathfrak{q}=\frac{2}{p}$ and, once the following inequality (used in the sequel without mentioning)
		\[
		(x+y)^s \leq 2^s (x^s+y^s) \quad \text{ for any } x,y \geq 0\;\textrm{ and }\; s>0
		\]
		is also considered, we obtain for every $t\in (0,\TM)$
		\begin{align}\label{a1}
	\frac{K^2}{2}	\left(\int_{\Omega} (u+1)^{(p+2\alpha-2)\theta}\right)^{\frac{1}{\theta}}&=\frac{K^2}{2} \|(u+1)^{\frac{p}{2}}\|_{L^{\frac{2(p+2\alpha-2)}{p}\theta}(\Omega)}^{\frac{2(p+2\alpha-2)}{p}}\\ \nonumber
		& \leq c_1 \|\nabla(u+1)^{\frac{p}{2}}\|_{L^2(\Omega)}^{\frac{2(p+2\alpha-2)}{p} a_1} \|(u+1)^{\frac{p}{2}}\|_{L^{\frac{2}{p}}(\Omega)}^{\frac{2(p+2\alpha-2)}{p} (1-a_1)}
		+ c_1 \|(u+1)^{\frac{p}{2}}\|_{L^{\frac{2}{p}}(\Omega)}^{\frac{2(p+2\alpha-2)}{p}},
		\end{align}
		where $c_1>0$ depends on $K$, $C_{GN}$, and with $a_1\in(0,1)$ taken from Lemma \ref{LemmaCoefficientAiAndExponents} and belonging to $(0,1)$ as therein proved. As a consequence, by observing that the mass conservation property \eqref{massConservation} implies the boundedness of $(u+1)^{\frac{p}{2}}$ in $L^{\infty}((0,T_{max}); L^{\frac{2}{p}}(\Omega))$,
		from \eqref{a1} we have that for some $c_2>0$  and $\beta_1\in(0,1)$ deduced from Lemma \ref{LemmaCoefficientAiAndExponents} 		
		\begin{equation} \label{Con_a1}
		\left(\int_{\Omega} (u+1)^{(p+2\alpha-2)\theta}\right)^{\frac{1}{\theta}}
		\leq c_2 \left(\int_{\Omega} |\nabla (u+1)^{\frac{p}{2}}|^2\right)^{\beta_1}+ c_2 \quad \textrm{ for every } \,t< T_{max}.
		\end{equation}
		In a similar way, we can again invoke the 
		Gagliardo--Nirenberg inequality, with an evident choice of $\mathfrak{p}$ and $\mathfrak{q}$, to have for some $c_3>0$ and $a_2\in (0,1)$ as in Lemma  \ref{LemmaCoefficientAiAndExponents}
		\begin{align*}
		&\left(\int_{\Omega} |\nabla v|^{2 \theta'}\right)^{\frac{1}{\theta'}} =\| |\nabla v|^q \|_{L^{\frac{2 \theta'}{q}}(\Omega)}^{\frac{2}{q}}\leq c_3 \|\nabla |\nabla v|^q\|_{L^2(\Omega)}^{\frac{2}{q}a_2} \| |\nabla v|^q\|_{L^{\frac{n}{q}}(\Omega)}^{\frac{2}{q}(1-a_2)}
		+ c_3 \| |\nabla v|^q\|_{L^{\frac{n}{q}}(\Omega)}^{\frac{2}{q}}\quad \textrm{ for all } t<T_{max}.
		\end{align*}
		In particular, by exploiting \eqref{Cg}, we entail that  (taking in mind $\gamma_1\in (0,1)$ from Lemma \ref{LemmaCoefficientAiAndExponents})		
		\begin{equation} \label{Con_a2}
		\left(\int_{\Omega} |\nabla v|^{2 \theta'}\right)^{\frac{1}{\theta'}} 
		\leq c_4 \left(\int_{\Omega} |\nabla |\nabla v|^q|^2 \right)^{\gamma_1}+ c_4\quad \textrm{ for all } t\in (0, T_{max}),
		\end{equation}
		with some computable $c_4>0$.
		
		Subsequently, by collecting \eqref{Y1}, \eqref{H1} and adjusting the product between \eqref{Con_a1} and \eqref{Con_a2} by means of the $\delta$-Young inequality, 	relation \eqref{Der.u} becomes
		\begin{align} \label{Der1.u}
		\frac{1}{p(p-1)}\frac{d}{dt} \int_{\Omega} (u+1)^p + \frac{1}{2} \int_{\Omega} (u+1)^{p-2} |\nabla u|^2 
		&\leq c_5 \left(\int_{\Omega} |\nabla (u+1)^{\frac{p}{2}}|^2\right)^{\beta_1}
		\left(\int_{\Omega} |\nabla |\nabla v|^q|^2 \right)^{\gamma_1}\\ \nonumber
		&\quad + \delta_1 \int_{\Omega} |\nabla (u+1)^{\frac{p}{2}}|^2 + \delta_2 \int_{\Omega} |\nabla |\nabla v|^q|^2 +c_5
		\end{align}
		for all $t\in(0,T_{max})$, arbitrary $\delta_1, \delta_2 >0$ and some $c_5>0$ depending also on $\delta_1$ and $\delta_2$.
		\subsubsection*{{\bf{Estimating $\frac{1}{q} \frac{d}{dt}\int_{\Omega} |\nabla v|^{2q}$ on $(0,T_{max})$.}}}
	First, by applying the identity $\Delta |\nabla v|^2 = 2 \nabla v \cdot \nabla \Delta v + 2 |D^2v|^2$,  we arrive for all $x \in \Omega$ and $t \in (0, T_{max})$ at 
		\[
		(|\nabla v|^2)_t = \Delta |\nabla v|^2 - 2 |D^2v|^2 - 2 |\nabla v|^2 + 2 \nabla g(u) \cdot \nabla v= 2 \nabla v \cdot \nabla \Delta v - 2 |\nabla v|^2 + 2 \nabla g(u) \cdot \nabla v.
		\]
		With such a relation in mind, by  using $|\nabla v|^{2q-2}$ as test function, a differentiation of the second equation of problem \eqref{problem}  implies that on $(0,T_{max})$ this estimate holds:
		\begin{align*}
		\frac{1}{q}\frac{d}{dt} \int_{\Omega} |\nabla v|^{2q} & = - (q-1) \int_{\Omega} |\nabla v|^{2q-4} |\nabla|\nabla v|^2|^2 
		+ \int_{\partial \Omega} |\nabla v|^{2q-2} \frac{\partial |\nabla v|^2}{\partial \nu} \\
		& \quad - 2  \int_{\Omega} |\nabla v|^{2q-2} |D^2 v|^2 - 2  \int_{\Omega} |\nabla v|^{2q} + 2 \int_{\Omega} |\nabla v|^{2q-2}  \nabla g(u) \cdot \nabla v.
		\end{align*}
		Now, an application of Lemma \ref{LemmaConv} allows us to find $C_{\eta} > 0$ such that for some suitable $\eta >0$ (to be chosen later) we have 
		\begin{align} \label{Der.v}
		&\frac{1}{q}\frac{d}{dt} \int_{\Omega} |\nabla v|^{2q}  + (q-1) \int_{\Omega} |\nabla v|^{2q-4} |\nabla|\nabla v|^2|^2 
		+ 2  \int_{\Omega} |\nabla v|^{2q-2} |D^2 v|^2 + 2  \int_{\Omega} |\nabla v|^{2q} \\ \nonumber
		&\quad \quad \quad \leq \eta \int_{\Omega} |\nabla |\nabla v|^q|^2 + C_{\eta} + 2 \int_{\Omega} |\nabla v|^{2q-2} \nabla g(u) \cdot \nabla v\quad \textrm{ on } \, (0,T_{max}).
		\end{align}
		By integrating by parts the latter integral above and using Young's inequality, we get
		\begin{align}\label{Yin}
		2 \int_{\Omega} |\nabla v|^{2q-2} \nabla g(u) \cdot \nabla v &= -2(q-1)\int_{\Omega} g(u) |\nabla v|^{2q-4} \nabla v \cdot \nabla|\nabla v|^2 
		- 2 \int_{\Omega} g(u)|\nabla v|^{2q-2} \Delta v \\ \nonumber
		& \leq \frac{(q-1)}{2} \int_{\Omega} |\nabla v|^{2q-4} |\nabla|\nabla v|^2|^2 + 2(q-1) \int_{\Omega}(g(u))^2 |\nabla v|^{2q-2}  \\ \nonumber
		&
		\quad + \frac{2}{n} \int_{\Omega} |\nabla v|^{2q-2} |\Delta v|^2 + \frac{n}{2} \int_{\Omega} (g(u))^2 |\nabla v|^{2q-2}  \textrm{ on }  (0,T_{max}),
		\end{align}
		where 
		\[
		\frac{2}{n} \int_{\Omega} |\nabla v|^{2q-2} |\Delta v|^2 \leq 2 \int_{\Omega} |\nabla v|^{2q-2} |D^2 v|^2
		\]
		in view of the pointwise inequality $|\Delta v|^2 \leq n |D^2 v|^2$. 
		Henceforth, by exploiting \eqref{Yin} and recalling assumption \eqref{AssumptionOnSegregationG}, we can rephrase \eqref{Der.v} as
		\begin{equation} \label{Der1.v}
		\frac{1}{q}\frac{d}{dt}\int_{\Omega} |\nabla v|^{2q} + \left(\frac{2(q-1)}{q^2} - \eta \right) \int_{\Omega} |\nabla |\nabla v|^q|^2 
		\leq c_6 \int_{\Omega} (u+1)^{2l} |\nabla v|^{2q-2}\quad  \textrm{ on }  (0,T_{max}),
		\end{equation}
		where $c_6$ is a positive constant depending also on $K_0$.  Let us now estimate the last integral in the previous bound. By  employing the H\"{o}lder inequality, we first obtain the following estimate
		\begin{equation} \label{H2}
		\int_{\Omega} (u+1)^{2l} |\nabla v|^{2q-2} \leq 
		\left(\int_{\Omega} (u+1)^{2l\mu}\right)^{\frac{1}{\mu}} \left(\int_{\Omega} |\nabla v|^{2(q-1)\mu'}\right)^{\frac{1}{\mu'}}\quad  \textrm{ on }  (0,T_{max}), 
		\end{equation}
		whereas by relying on Lemma \ref{LemmaGN}, we find a constant $c_7>0$, depending  on $C_{GN}$, such that for $a_3 \in (0,1)$ from Lemma  \ref{LemmaCoefficientAiAndExponents} we arrive at
		\begin{align}\label{a3}
		\left(\int_{\Omega} (u+1)^{2l\mu}\right)^{\frac{1}{\mu}}&= \|(u+1)^{\frac{p}{2}}\|_{L^{\frac{4\mu l}{p}}(\Omega)}^{\frac{4l}{p}}\\ \nonumber
		& \leq c_7 \|\nabla(u+1)^{\frac{p}{2}}\|_{L^2(\Omega)}^{\frac{4l}{p}  a_3} \|(u+1)^{\frac{p}{2}}\|_{L^{\frac{2l}{p}}(\Omega)}^{\frac{4l}{p} (1-a_3)}
		+ c_7 \|(u+1)^{\frac{p}{2}}\|_{L^{\frac{2l}{p}}(\Omega)}^{\frac{4l}{p}}\quad \textrm{on } \, (0,\TM).
		\end{align}
		On the other hand, by arguing as before, we infer for some $c_8>0$, $\beta_2\in (0,1)$ taken by Lemma \ref{LemmaCoefficientAiAndExponents}
		and by the finiteness  of $\|(u+1)^{\frac{p}{2}}\|_{L^{\frac{2l}{p}}(\Omega)}$ (immediately coming from \eqref{massConservation} in view of $0<l<\frac{2}{n}<1$) 
		\begin{equation} \label{Con_a3}
	\left(\int_{\Omega} (u+1)^{2l\mu}\right)^{\frac{1}{\mu}}
	\leq c_8 \left(\int_{\Omega} |\nabla (u+1)^{\frac{p}{2}}|^2\right)^{\beta_2}+ c_8\quad \textrm{for all }\, t \in (0,T_{max}).
	\end{equation}
		(Let us note that in \eqref{a3} we have intentionally applied the Gagliardo--Nirenberg inequality with exponent $\mathfrak{q}=\frac{2l }{p}$ only for exhibiting reasons; in particular, the expression of $\mathcal{B}$ in Lemma \ref{pqtildeLemman=2andn>2Lemma} appears more compact than the one that would be obtained by considering the optimal exponent $\mathfrak{q}=\frac{2}{p}$. This does not preclude the sharpness of the assumption because, since $\mu$ is  taken indefinitely large, the exponent $\mathfrak{p}$ has the control on $a_3$, and not $\mathfrak{q}$.)
		
		At this point, by making use again of the Lemma \ref{LemmaGN}, positive constants $c_9$ and $c_{10}$ satisfy 
		\begin{align*}
		\left(\int_{\Omega} |\nabla v|^{2(q-1) \mu'}\right)^{\frac{1}{\mu'}} &=\| |\nabla v|^q \|_{L^{\frac{2(q-1)}{q}\mu'}(\Omega)}^{\frac{2(q-1)}{q}}\\ 
		& \leq c_9 \|\nabla |\nabla v|^q\|_{L^2(\Omega)}^{\frac{2(q-1)}{q}  a_4} \| |\nabla v|^q\|_{L^{\frac{n}{q}}(\Omega)}^{\frac{2(q-1)}{q}  (1-a_4)}
		+ c_9 \| |\nabla v|^q\|_{L^{\frac{n}{q}}(\Omega)}^{\frac{2(q-1)}{q}}\quad \textrm{ on } (0,T_{max}),
		\end{align*}
	\begin{equation} \label{Con_a4}
	\left(\int_{\Omega} |\nabla v|^{2(q-1)\mu'}\right)^{\frac{1}{\mu'}} 
	\leq c_{10} \left(\int_{\Omega} |\nabla |\nabla v|^q|^2 \right)^{\gamma_2}+ c_{10}\quad \textrm{ for every } t\in (0,T_{max}),
	\end{equation}
	where, once more through Lemma \ref{LemmaCoefficientAiAndExponents},  $\gamma_2\in (0,1)$ and $a_4 \in (0,1).$
		
		Finally, by plugging relations \eqref{H2}, \eqref{Con_a3} and \eqref{Con_a4} into  bound
		\eqref{Der1.v}, a further application of the  $\delta$-Young inequality
	\begin{align} \label{Der2.v}
	\frac{1}{q}\frac{d}{dt}\int_{\Omega} |\nabla v|^{2q} + \left(\frac{2(q-1)}{q^2} - \eta \right) \int_{\Omega} |\nabla |\nabla v|^q|^2 & \leq c_{11} \left(\int_{\Omega} |\nabla (u+1)^{\frac{p}{2}}|^2\right)^{\beta_2}
	\left(\int_{\Omega} |\nabla |\nabla v|^q|^2 \right)^{\gamma_2}\\ \nonumber
	&\quad + \delta_3 \int_{\Omega} |\nabla (u+1)^{\frac{p}{2}}|^2 + \delta_4 \int_{\Omega} |\nabla |\nabla v|^q|^2 
	+ c_{11}\quad \textrm {for all }\, t< T_{max},
	\end{align}
	with $\delta_3, \delta_4>0$ and $c_{11}$ being a proper positive constant which depends on $\delta_3$ and $\delta_4$.
		\subsubsection*{{\bf{Combining terms: the absorptive inequality on $(0,T_{max})$}}}
		By adding the two contributions from \eqref{Der1.u} and \eqref{Der2.v}, yields for some $c_{12}>0$ 
			\begin{align} \label{Der}
	&	\frac{d}{dt}\left(\frac{1}{p(p-1)}\int_{\Omega} (u+1)^p + \frac{1}{q} \int_{\Omega} |\nabla v|^{2q}\right) + \left(\frac{2}{p^2}-(\delta_1+\delta_3)\right) 
		\int_{\Omega} |\nabla(u+1)^{\frac{p}{2}}|^2 \\ \nonumber 
		& +
		\left(\frac{2(q-1)}{q^2} -(\eta + \delta_2+\delta_4)\right) \int_{\Omega} |\nabla|\nabla v|^q|^2 \leq c_{12} \left(\int_{\Omega} |\nabla (u+1)^{\frac{p}{2}}|^2\right)^{\beta_1} \left(\int_{\Omega} |\nabla |\nabla v|^q|^2 \right)^{\gamma_1}
		\\ \nonumber
		&+ c_{12} \left(\int_{\Omega} |\nabla (u+1)^{\frac{p}{2}}|^2\right)^{\beta_2} \left(\int_{\Omega} |\nabla |\nabla v|^q|^2 \right)^{\gamma_2} + c_{12}\quad \textrm{ on }  (0,\TM),
		\end{align}
		where accordingly to Lemma \ref{LemmaCoefficientAiAndExponents}, the coefficients 
		$\beta_1+\gamma_1\in (0,1)$ and $\beta_2+\gamma_2\in (0,1)$.
		% fulfill
		%\[
		%\beta_1 + \gamma_1 = 	\frac{p-2+2\alpha}{p}a_1+\frac{1}{q}a_2 \in (0,1)
		%\,\textrm{ and }\,
		%\
		%\beta_2+\gamma_2 = \frac{2 l }{p}a_3+\frac{q-1}{q}a_4   \in (0,1).
		%\]
		Therefore we can apply the first inequality of Lemma \ref{LemmaIneq} to \eqref{Der},  where we choose $\delta_1=\delta_3=\frac{1}{2p^2}$, 
		$\delta_2=\delta_4=\eta=\frac{(q-1)}{3q^2}$, so to find a positive constant $c_{13}$ producing
		\begin{align} \label{Der1}
		&\frac{d}{dt}\left(\frac{1}{p(p-1)}\int_{\Omega} (u+1)^p + \frac{1}{q} \int_{\Omega} |\nabla v|^{2q}\right) + \frac{1}{p^2} \int_{\Omega} |\nabla(u+1)^{\frac{p}{2}}|^2 +\frac{(q-1)}{q^2} \int_{\Omega} |\nabla|\nabla v|^q|^2\leq c_{13}\quad \textrm{on }  (0,T_{max}).
		%
		%  \\ \nonumber
		%&\leq \frac{1}{p^2} \int_{\Omega} |\nabla(u+1)^{\frac{p}{2}}|^2 + \left[\frac{(q-1)}{2q^2}\right] \int_{\Omega} |\nabla|\nabla v|^q|^2 + 
		\end{align}
		Again by employing twice the Gagliardo--Nirenberg inequality, we have for $\kappa_1 \in(0,1)$ and $\kappa_2 \in(0,1)$ derived in Lemma \ref{LemmaCoefficientAiAndExponents}, and suitable large $c_{14}>0$, that these estimates hold true for all $t\in(0,T_{max})$:
		\[
		\int_{\Omega} (u+1)^p= \|(u+1)^{\frac{p}{2}}\|_{L^2(\Omega)}^2  \leq c_{14} \|\nabla(u+1)^{\frac{p}{2}}\|_{L^2(\Omega)}^{2 \kappa_1}  
		\|(u+1)^{\frac{p}{2}}\|_{L^{\frac{2}{p}}(\Omega)}^{2 (1-\kappa_1)} + c_{14} \|(u+1)^{\frac{p}{2}}\|_{L^{\frac{2}{p}}(\Omega)}^2,
		\]
		and
		\[
		\int_{\Omega} |\nabla v|^{2q} =\| |\nabla v|^q \|_{L^2(\Omega)}^2 \leq c_{14} \|\nabla |\nabla v|^q\|_{L^2(\Omega)}^{2 \kappa_2} 
		\| |\nabla v|^q\|_{L^{\frac{n}{q}}(\Omega)}^{2 (1-\kappa_2)} + c_{14} \| |\nabla v|^q\|_{L^{\frac{n}{q}}(\Omega)}^2.
		\]
		The already used mass conservation property and the boundedness of $\lVert v(\cdot,t)\rVert_{W^{1,n}(\Omega)}$, provide some positive constant $c_{15}$ such that 
		\begin{equation}\label{K1}
		\int_{\Omega} (u+1)^p \leq c_{15} \left(\int_{\Omega} |\nabla(u+1)^{\frac{p}{2}}|^2\right)^{\kappa_1} + c_{15} \quad \textrm{ for all } t\in (0,T_{max}),
		\end{equation}
		and
		\begin{equation}\label{K2}
		\int_{\Omega} |\nabla v|^{2q} \leq c_{15} \left(\int_{\Omega} |\nabla |\nabla v|^q|^2\right)^{\kappa_2}+ c_{15} \quad \textrm{ on }\, (0,T_{max}).
		\end{equation}
		Consequently, by collecting \eqref{K1} and \eqref{K2}, we can rewrite \eqref{Der1} in the following way
		\[
		\frac{d}{dt} \left(\frac{1}{p(p-1)}\int_{\Omega} (u+1)^p + \frac{1}{q} \int_{\Omega} |\nabla v|^{2q} \right) + c_{16} \left(\int_{\Omega} (u+1)^p\right)^{\frac{1}{\kappa_1}} 
		+ c_{16} \left(\int_{\Omega} |\nabla v|^{2q}\right)^{\frac{1}{\kappa_2}} \leq c_{17}\quad \textrm{ on }\,  (0,T_{max}),
		\]
		with  positive constants $c_{16}, c_{17}$.
		
		From all of the above, we invoke the second inequality in Lemma \ref{LemmaIneq}, so to see that the function $y=y(t)$ satisfies this initial value problem 
		\[
		\begin{cases}
		y'(t) + c_{18} y^{\kappa}(t) \leq c_{19}\quad \textrm{ on } (0,T_{max}), \\
		y(0)=y_0=\frac{1}{p(p-1)}\int_{\Omega} (u_0(x)+1)^p dx + \frac{1}{q} \int_{\Omega} |\nabla v_0(x)|^{2q} dx,
		\end{cases}
		\]
		with suitable constants $\kappa, c_{18}, c_{19}>0$. This leads to the conclusion for appropriate $L>0$ since standard ODE comparison arguments give
		\[
		y(t) \leq \max\left\{y_0, \left(\frac{c_{19}}{c_{18}}\right)^{\frac{1}{\kappa}}\right\}\quad \textrm{ for every } \, t<T_{max}.
		\]
	\end{proof}
\end{lemma}
With these gained bounds, we exploit a  general boundedness result to quasilinear parabolic equations (see \cite{TaoWinkParaPara}) so to ensure uniform-in-time boundedness of the local solution $(u,v)$ to system \eqref{problem}.
%\section{Proof of the claims}\label{SectionProofsClaims}
\subsubsection*{{\bf{Proof of Theorem \ref{MainTheorem}}}} Let $(u,v)$ be the local classical solution to \eqref{problem}. Upon enlarging $p$ and $q$ accordingly to what said in item \ref{ItemUponEnlarging}  of Remark \ref{RemarkEnlargingPandQ}, we can obtain that the term $f(u)\nabla v\in L^\infty((0,\TM);L^{q_1}(\Omega))$, for some $q_1>n+2$.  So we conclude thanks to Lemma \ref{LemmaStime}, \cite[Lemma A.1]{TaoWinkParaPara}  and the dichotomy criterion \eqref{ExtensibilityCrit}.
\qed
\begin{remark}[\rm{Some hints about the parabolic-elliptic model}]\label{HintsEllipticCase} 
Let us consider the equations
	\begin{equation}\label{problemSimlified}
	u_t= \Delta u - \nabla \cdot (f(u) \nabla v) \; \textrm{ and }\;
0=\Delta v-v+g(u),  \quad  \text{ in } \Omega \times (0,T_{max}),
\end{equation}
endowed with homogeneous Neumann boundary conditions, nontrivial initial data $u(x,0)=u_0(x)\geq 0$, where $f$ and $g$ comply with assumptions in Theorem \ref{MainTheorem}. Similarly to what already done, we have
\begin{align*}%\label{Der.uSimplified}
\frac{1}{p(p-1)}\frac{d}{dt} \int_{\Omega} (u+1)^p &= \frac{1}{p-1}\int_{\Omega} (u+1)^{p-1} \nabla \cdot \nabla u - \frac{1}{p-1}
\int_{\Omega} (u+1)^{p-1} \nabla \cdot (f(u) \nabla v) \\ \nonumber
%&= - \int_{\Omega} \phi''(u) |\nabla u|^2 + \int_{\Omega} \phi''(u) f(u) \nabla u \cdot \nabla v \\ \nonumber
&= - \int_{\Omega} (u+1)^{p-2} |\nabla u|^2 + \int_{\Omega} (u+1)^{p-2} f(u) \nabla u \cdot \nabla v  \\ \nonumber 
& \leq - \int_{\Omega} (u+1)^{p-2} |\nabla u|^2 + \frac{K K_0}{p+\alpha-1}\int_{\Omega} (u+1)^{p+\alpha+l-1} \quad \textrm{ on }\, (0, T_{max}).
\end{align*}
This estimate is essentially the same than that derived in \cite[$\S$4]{WinklerNoNLinearanalysisSublinearProduction} so that, as therein, in order to take advantage from a combination of the Gagliardo--Nirenberg and Young's inequalities, one has to impose $\alpha-1+l<\frac{2}{n}$ (coinciding, exactly as discussed in $\S$\ref{SubSectionOvervieParaEll}, with the parabolic-elliptic version of assumption \eqref{AssumptionAlphaElleN}); consequently, the integral $\frac{K K_0}{p+\alpha-1}\int_{\Omega} (u+1)^{p+\alpha+l-1}$ can be suitably treated. Standard procedures, successively, provide that $u\in L^\infty((0,\TM);L^p(\Omega))$ for arbitrarily large $p>1$, and hence also $g\in L^\infty((0,\TM);L^p(\Omega))$ for any $l\in(0,1).$ Finally, elliptic regularity theory applied to the second equation in problem \eqref{problemSimlified}  infers uniform bound of $\nabla v$, on $(0,\TM)$, so that $u$ and $v$ are uniformly bounded for all $t>0.$ 

We note that the necessary regularity of $\nabla v$ is gained only by differentiating $\int_{\Omega} (u+1)^p$, by using the initial-boundary value problem \eqref{problemSimlified} and, solely, the mass conservation property; neither an estimate like that in \eqref{Cg}  is a priory needed nor the analysis of the term $\int_\Omega |\nabla v|^{2q}$, involving the extra parameter $q$, has to be developed.
\end{remark}
\subsubsection*{Acknowledgments}
%The author is grateful to the referees for helpful suggestions which improved this article.
The authors are members of the Gruppo Nazionale per l'Analisi Matematica, la Probabilit\`a e le loro Applicazioni (GNAMPA) of the Istituto Na\-zio\-na\-le di Alta Matematica (INdAM) and  are  partially supported by the research project \textit{Integro-differential Equations and Non-Local Problems}, funded by Fondazione di Sardegna (2017). GV is also partially supported by MIUR (Italian Ministry of Education, University and Research) Prin 2017 \textit{Nonlinear Differential Problems via Variational, Topological and Set-valued Methods} (Grant Number: 2017AYM8XW). 

% % % % % % %

%\bibliography{Bibliography}{}
%\bibliography{reference}
%\bibliographystyle{abbrv}
%\bibliographystyle{unsrt}
%\bibliographystyle{acm}
%\bibliographystyle{apalike}
%\bibliographystyle{alpha}
%\bibliographystyle{AIMS} 
%\bibliographystyle{plainnat}
\end{document}